\documentclass[12pt,reqno]{amsart}
\usepackage{hyperref}
\usepackage{amsmath,amssymb,physics}
\usepackage{amsmath,amsthm,amsopn,amssymb,a4wide,varioref}
\usepackage{enumitem}
\usepackage[mathscr]{eucal}
\usepackage{color, bm, amscd, tikz-cd}
\usepackage{nicematrix} 
\newcommand{\ran}{\operatorname{ran}}
\newcommand{\spn}{\operatorname{span}}
\newcommand{\spncl}{\overline{\operatorname{span}}}

\newcommand{\C}{\mathbb{C}}
\newcommand{\Z}{\mathbb{Z}}

\newcommand{\D}{{\mathbb D}}
\newcommand{\Dc}{\overline{\mathbb D}}
\newcommand{\T}{\mathbb{T}}
\newcommand{\B}{\mathscr{B}}
\renewcommand{\H}{\mathcal{H}}
\newcommand{\K}{\mathcal{K}}
\newcommand{\E}{\mathcal{E}} 
\newcommand{\F}{\mathcal{F}} 
\renewcommand{\L}{\mathcal{L}} 
\newcommand{\W}{\mathcal W}
\newcommand{\M}{\mathcal{M}}

\newcommand{\la}{ \langle }
\newcommand{\ra}{\rangle}

\newtheorem{thm}{Theorem}[section]

\newtheorem{corollary}[thm]{Corollary}
\newtheorem{lemma}[thm]{Lemma}

\newtheorem{proposition}[thm]{Proposition}

\newtheorem{definition}[thm]{Definition}
\newtheorem{remark}[thm]{Remark}
\newtheorem{note}[thm]{Note}
\newtheorem{example}[thm]{Example}

\numberwithin{equation}{section}

\def\textmatrix#1&#2\\#3&#4\\{\bigl({#1 \atop #3}\ {#2 \atop #4}\bigr)}
\def\dispmatrix#1&#2\\#3&#4\\{\left({#1 \atop #3}\ {#2 \atop #4}\right)}
\numberwithin{equation}{section}

\def\textmatrix#1&#2\\#3&#4\\{\bigl({#1 \atop #3}\ {#2 \atop #4}\bigr)}
\def\dispmatrix#1&#2\\#3&#4\\{\left({#1 \atop #3}\ {#2 \atop #4}\right)}

\begin{document}
	
	\title[A pair of commuting isometries]{The joint spectrum for a commuting pair of  isometries in certain cases}

	\author{Tirthankar Bhattacharyya}
	\address{Department of Mathematics, Indian Institute of Science, Bangalore 560012.}
	\email{tirtha@iisc.ac.in}
	
	\author{Shubham Rastogi}
	\address{Department of Mathematics, Indian Institute of Science, Bangalore 560012.}
	\email{shubhamr@iisc.ac.in}
	
	\author{Vijaya Kumar U.}
	\address{Department of Mathematics, Indian Institute of Science, Bangalore 560012.}
	\email{vijayak@iisc.ac.in}
	
	\maketitle

\vspace*{2mm}

\begin{centerline}
{In fond memory of J\"{o}rg Eschmeier}
\end{centerline}
	
	\renewcommand{\thefootnote}{\fnsymbol{footnote}}
	
	
	\begin{abstract}
	We show that the joint spectrum of two  commuting isometries can vary widely depending on various factors. It can range from being small (of measure zero or an analytic disc for example) to the full bidisc. En route, we discover a new model pair in the negative defect case and relate it to the modified bi-shift.
	\end{abstract}
	
	\section{Introduction}
	An isometry $V$ is called $pure$ if $V^{*n}$ converges to $0$ strongly as $n \rightarrow \infty$. This is equivalent to saying that $V$ is the unilateral shift of multiplicity equal to the dimension of the range of the {\em defect operator} $I - VV^*$.
	
	The famous Wold decomposition \cite{von Neumann, Wold} tells us that given an isometry $V$ on a Hilbert space $\mathcal H$,  the space $\mathcal H$ breaks uniquely into a direct sum $\mathcal H = \mathcal H_0 \oplus \H_0^\perp$ of reducing subspaces such that $V|_{\H_0}$ is a unitary and $V|_{\H_0^\perp}$ is a pure isometry. This immediately implies that for a non-unitary isometry $V$ (i.e., when the defect operator is positive and not zero), the spectrum $\sigma(V)$ is the  closed unit disc $\Dc=\{z\in \C:|z|\le 1\}.$ The situation for a pair of commuting isometries is vastly different.
	
The topic of commuting isometries has been vigorously pursued in the last two decades, see \cite{A-K-M, B-D-F, Burdak, BKPS, BKS,GG, Kosiek, Sarkar, Popovici I, Popovici, Slocinski} and the references therein.  In \cite{BKS} and \cite{BKPS}, the novel idea of using graphs has led to a clear understanding of structures. 


 The	{\em defect operator} $ C(V_1,V_2)$ is introduced in \cite{Guo} and \cite{Yang},  as
	\begin{equation*}
		C(V_1,V_2)=I-V_1V_1^*-V_2V_2^*+V_1V_2V_2^*V_1^*.
\end{equation*}
	In \cite{Yang} and \cite{Sarkar}, the authors provide  the characterization of $(V_1,V_2),$ when the defect is positive, negative or zero.  It is well known (see \cite{Yang,GG}) that a pair has positive defect if and only if it is doubly commuting, and it has negative  defect if and only if it is  dual doubly commuting.
	 In all the three cases, the defect is either a projection or negative of a projection. In general the defect is the difference of two projections; see \cite{Sarkar}.

	In this paper we study the pairs of commuting isometries, whose defect is the difference of two mutually orthogonal projections. We  characterize such pairs in Theorem \ref{thm:iff} and we classify them in Table \ref{classification:table}. We also provide the characterization for a few cases in Table \ref{classification:table}; see Lemma \ref{char off} and Lemma \ref{lem:rangeincl}.  We rephrase the structure of $(V_1,V_2)$ in each case, which appears in Table \ref{classification:table}, in a unified approach using the Berger-Coburn-Lebow (BCL) Theorem. The joint spectrum is studied in detail for all the cases except the last one appearing in Table \ref{classification:table}.

	There is the related concept of the {\em fringe operators}:
	$$F_1: \ker V_1^*\rightarrow \ker V_1^* \text{ and } F_2: \ker V_2^*\rightarrow \ker V_2^*$$
	defined  by,
	\begin{equation}
		F_1(x)=P_{\ker V_1^*}V_2(x)\text{ and } F_2(x)=P_{\ker V_2^*}V_1(x).
	\end{equation}
	In various characterizations of Table \ref{classification:table}, we shall point out the criteria in terms of the fringe operators for possible use in examples.

	\subsection{The joint spectrum}
	If $(T_1,T_2)$ is a pair of commuting bounded operators on $\H,$ then for defining (see \cite{Kallstrom, Taylor}) the {\em Taylor joint spectrum} $\sigma(T_1,T_2),$ one considers the {\em Koszul complex} $K(T_1,T_2):$
	\begin{equation}\label{eq:Kozul}
		0\overset{\delta_0}{\to} \H \overset{\delta_1}{\to}\H\oplus \H\overset{\delta_2}{\to}\H\overset{\delta_3}{\to} 0
	\end{equation}
	where $\delta_1(h)= (T_1h,T_2h)$ for $h\in \H$ and $\delta_2(h_1,h_2)= T_1h_2-T_2h_1$ for $h_1,h_2\in \H.$	From the way the complex is constructed, $\ran \delta_{n-1} \subseteq \ker \delta_{n}.$ When $\ran \delta_{n-1}= \ker \delta_{n}$ for all $n=1,2,3$ we say that the {\em Koszul complex $K(T_1,T_2)$ is exact or the pair $(T_1,T_2)$ is non-singular.}
	A pair $(\lambda_1,\lambda_2)\in \C^2$ is said to be in the {\em joint spectrum} $\sigma(T_1,T_2)$ if the pair $(T_1-\lambda_1I, T_2-\lambda_2I)$ is singular. In the case of a singular pair, we say that the {\em non-singularity breaks at the stage $n$} if $\ran \delta_{n-1} \neq \ker \delta_{n}$.
	
	Observe that the non-singularity breaks at stage 1 if and only if $(\lambda_1,\lambda_2)$ is a joint eigenvalue for $(T_1,T_2)$ and the non-singularity breaks at stage 3 if and only if the joint range of $(T_1-\lambda_1I, T_2-\lambda_2I)$ is not the whole space $\H.$ If $(\overline{\lambda_1},\overline{\lambda_2})$ is a joint eigenvalue of $(T_1^*,T_2^*)$, then by the fact that $\ran T_1 + \ran T_2 = \H$ implies  $\ker T_1^* \cap \ker T_2^* = \{0\}$, the non-singularity of the Koszul complex $K(T_1-\lambda_1I,T_2-\lambda_2I)$ breaks at stage 3. There are a few elementary results which we record as a lemma so that we can refer to it later.
	\begin{lemma} \label{ReducingSpectrum}
		Let $\H$ and $\K$ be two non-zero Hilbert spaces. Let $(T_1,T_2)$ be a pair of commuting bounded operators on $\H.$
		\begin{enumerate}
			\item $\sigma(T_1, T_2) \subseteq \sigma(T_1) \times \sigma(T_2)$.\label{Part1ReducingSpectrum}
			
			\item\label{Part2ReducingSpectrum}
			If there is a non-trivial joint reducing subspace $\H_0$  for $(T_1,T_2),$
			i.e., if $\H=\H_0\oplus\H_0^\perp$ and
			\[T_i=\begin{pNiceMatrix}[first-row,last-col,nullify-dots]
				\H_0& \H_0^\perp\\
				T_{i0} & 0 & \H_0\\
				0& T_{i1} & \H_0^\perp
			\end{pNiceMatrix}\] 
			then
			$$\sigma(T_1,T_2)=\sigma(T_{10},T_{20})\cup \sigma(T_{11},T_{21}).$$
			
			\item $(z_1, z_2) \in \sigma(T_1, T_2)$ if and only if $(\overline{z_1}, \overline{z_2}) \in \sigma (T_1^*, T_2^*)$.
			
			\item  $\sigma(I_\K\otimes T_1,I_\K\otimes T_2)=\sigma(T_1, T_2)=\sigma(T_1\otimes I_\K,T_2\otimes I_\K).$\label{Part4ReducingSpectrum}
			
			\item Let $(S_1,S_2)$ be a pair of commuting bounded operators on a Hilbert space $\K.$ If $(T_1,T_2)$ is jointly unitarily equivalent to $(S_1,S_2),$ then $\sigma(T_1,T_2)=\sigma (S_1,S_2).$
			
			\item For any $T$ in $\B(\H)$ and $S$ in $\B(\K)$, the joint spectrum $\sigma(T \otimes I_\K , I_\H \otimes S)$ is the Cartesian product $\sigma(T) \times \sigma(S)$.\label{Part6ReducingSpectrum}

		\end{enumerate}
	\end{lemma}

	Thus, for commuting isometries $V_1$ and $V_2$, we have $\sigma(V_1, V_2) \subseteq \overline{\D^2}$. The joint spectrum of a pair of commuting unitary operators is contained in the torus $\T^2,$ where $\T=\{z\in\C: |z|=1\}$ is the unit circle in the complex plane.

		This note uses the {\em fundamental pairs of isometries consisting of multiplication operators} to describe the structure of $(V_1,V_2).$ These are fundamental in the sense that in each case the sign of the defect operator is dictated by the fundamental pair alone.

	When the defect operator $C(V_1,V_2)$ of two commuting isometries is positive or negative, but not zero, then the whole space $\H$ breaks into a direct sum of reducing subspaces $\H = \H_0 \oplus \H_0^\perp$ in the style of Wold where the restriction of $(V_1,V_2)$ on the $\H_0$ part  is the fundamental pair and the restriction of $(V_1,V_2)$ on the $\H_0^\perp$ part has defect zero. What the fundamental pair is depends on whether $C(V_1,V_2)$ is positive or negative. These are the contents of Theorem \ref{StructureNegative} and Theorem \ref{StructurePositive}.
	
	The fundamental pairs are such that in both cases (of $C(V_1,V_2)$  positive or negative), the joint spectrum of $(V_1,V_2)$ is the whole closed bidisc $\overline{\mathbb D^2}$. These are done in Theorem \ref{StructureNegative} and Theorem \ref{JSPositive}.	If the defect operator $C(V_1, V_2)$ is zero,  the joint spectrum of $(V_1,V_2)$ is contained in the topological boundary of the bidisc. 
	
	The structure theorem in the case $\ran V_1=\ran V_2$, shows that $(V_1, V_2)$ is the direct sum of a prototypical pair (see Subsection \ref{ProtypicalOffDiag}) and a pair of commuting unitaries. The joint spectrum is computed. The joint spectrum of the prototypical pair of this case, is  neither the closed bidisc nor contained inside the topological boundary of the bidisc.
	
  In the  case
  $\ran V_2\subsetneq\ran V_1,$	the joint spectrum
$\sigma(V_1,V_2)\subseteq \{(z_1,z_1z_2): z_1,z_2\in \Dc\}.$ The above inclusion is sharp; see Example \ref{eg:sharp}, and it can be a strict inclusion; see Example \ref{eg:strict}. Note that $\{(z_1,z_1z_2): z_1,z_2\in \Dc\}$ has measure non-zero and it is not equal to the closed bidisc.

	In each case above except the case $\ran V_2\subsetneq \ran V_1,$ we point out the stage of the Koszul complex where non-singularity is broken. 
	
	\subsection{The Berger-Coburn-Lebow Theorem}
	For a Hilbert space $\mathcal E$, the Hardy space of $\mathcal E$-valued functions on the unit disc in the complex plane is
	$$H^2_{\mathbb D}(\mathcal E) = \{ f : \mathbb D \rightarrow \mathcal E \mid f \text{ is analytic and } f(z) = \sum_{n=0}^{\infty} a_n z^n \text{ with } \sum_{n=0}^{\infty} \|a_n\|_\E^2 < \infty\}.$$
	Here the $a_n$ are from $\mathcal E$. This is a Hilbert space with the inner product
	$$\langle \sum_{n=0}^{\infty} a_n z^n , \sum_{n=0}^{\infty} b_n z^n \rangle =  \sum_{n=0}^{\infty} \langle a_n , b_n \rangle_\E$$ and is identifiable with $H^2_{\mathbb D} \otimes \mathcal E$ where $H^2_{\D}$ stands for the Hardy space of scalar-valued functions on $\mathbb D$. We shall use this identification throughout the paper, often without any further mention, and $M_z$ denotes the multiplication operator by the coordinate function $z$ on $H^2_{\mathbb D}.$    For $\lambda \in \D$, let $k_\lambda$ be the function in $H^2_{\D}$ given by
	$$k_\lambda(z) = \sum_{n=0}^\infty z^n \overline{\lambda}^n = \frac{1}{1 - z\overline{\lambda}}.$$
	It is well-known that the span of $\{ k_\lambda : \lambda \in \D\}$ is dense in $H^2_{\D}$.
	
	The space of $\B(\E)$-valued bounded analytic functions on $\D$ will be denoted by $H^\infty_\D(\B(\E)).$ 	Naturally, if $\varphi\in H^\infty_\D(\B(\E))$, then it induces a multiplication operator $M_\varphi$ on $H^2_\D(\E)$.  One of the main tools for us is the Berger-Coburn-Lebow (BCL) theorem \cite{Berger}.
	
	\begin{thm}\label{thm:BCL}
		Let $(V_1,V_2)$ be a commuting pair of  isometries acting  on $\H.$ Then, up to unitary equivalence, the Hilbert space $\H$ breaks into a direct sum of reducing subspaces $\H=\H_p\oplus \H_u$ such that
		\begin{enumerate}
			\item There is a unique (up to unitary equivalence) triple $(\E, P, U)$ where  $\E$ is a Hilbert space, $P$ is a projection on $\E$ and $U$ is a unitary on $\E$ such that $\H_p=H^2_{\D}(\E),$ the functions $\varphi_1$ and $\varphi_2$ defined on $\D$ by
			\begin{equation}\label{eq:mult}
				\varphi_1(z)=U^*(P^\perp+zP) \text{ and } \varphi_2(z)=(P+zP^\perp)U,
			\end{equation}
			are commuting multipliers in $H^\infty_\D(\B(\E))$ and $(V_1|_{\H_p},V_2|_{\H_p})$ is equal to $(M_{\varphi_1},M_{\varphi_2}).$
			
			\item    $V_1|_{\H_u}$ and $V_2|_{\H_u}$ are commuting unitary operators.
		\end{enumerate}
	\end{thm}
	
	The result of the theorem above will be called the {\em BCL representation} of $(V_1, V_2)$. Using Theorem \ref{thm:BCL} we can compute the defect operator (see \cite{Yang}), because  $M_{\varphi_1}=I_{H^2_{\D}}\otimes U^*P^{\perp}+ M_z \otimes U^*P$ and $M_{\varphi_2}=I_{H^2_{\D}}\otimes PU+ M_z \otimes P^\perp U.$ Hence	$$C(M_{\varphi_1}, M_{\varphi_2})= (I- M_zM_z^*)\otimes (U^*PU-P)= E_0 \otimes (U^*PU-P)$$
	where $E_0$ is the {\em one dimensional projection onto the space of constant functions} in $H^2_\D$. Together with the fact that the defect operator of a pair of commuting unitary operators is zero, this means, in the decomposition $\H=\H_p\oplus \H_u$ with $\H_p=H^2_{\D}(\E),$
	\begin{equation} \label{DefectandBCL} C(V_1, V_2) = (E_0\otimes (U^*PU-P))\oplus 0. \end{equation}

	\begin{definition}
		
		\leavevmode
		
		\begin{enumerate}
			\item A {\em BCL triple} $(\E,P,U)$ is a Hilbert space $\E,$ along with a projection $P$ and a unitary $U.$ It is said to be the {\em BCL triple for the pair of commuting isometries $(V_1,V_2)$} if $\E,P$ and $U$ are as in Theorem \ref{thm:BCL}, part  (1).
			
			\item Given a BCL triple $(\E,P,U),$  the
			functions $\varphi_1,\varphi_2:\D\to \B(\E)$ will always be as defined in \eqref{eq:mult}.
			
			\item A pair of commuting isometries $(V_1,V_2)$ is called pure if $\mathcal H_u = \{0\}$ in its BCL representation.
		\end{enumerate}
	\end{definition}
	
	In fact, $\H_u$ is the unitary part of the product $V=V_1V_2$ in its Wold decomposition and $\E=\ker V^*$. Thus  $\H_u=\{0\}$ if and only if $V$ is pure.
	
	The Berger-Coburn-Lebow theorem has an interesting consequence in the case when the part $\mathcal H_p$ is non-zero. We have $\varphi_1(z)\varphi_2(z) = z$ for every $z \in \mathbb D$ and hence $M_{\varphi_1} M_{\varphi_2} = M_z\otimes I_\E$. Hence, by the spectral mapping theorem for joint spectra (see for example \cite{Curto}),
	$$ \{ z_1z_2 : (z_1, z_2) \in \sigma (M_{\varphi_1}, M_{\varphi_2}) \} = \sigma(M_z\otimes I_\E) = \overline{\mathbb D}.$$
	Thus, if $\mathcal H_p \neq \{0\}$, then $\sigma (M_{\varphi_1}, M_{\varphi_2})$ cannot be contained in the torus $\T^2$. Hence by Lemma \ref{ReducingSpectrum} part \eqref{Part2ReducingSpectrum}, $\sigma (V_1,V_2)$ cannot be contained in the torus $\T^2.$

	Let $(\E,P,U)$ be a BCL triple. It is easy to see that if  the Koszul complex $K(\varphi_1(z)-\lambda_1I,\varphi_2(z)-\lambda_2I)$ breaks at stage 3 for some $z\in \D,$ then the Koszul complex $K(M_{\varphi_1}-\lambda_1I,M_{\varphi_2}-\lambda_2I)$ breaks at stage 3. More generally, we show that, in all the cases under consideration in this note except the case $\ran V_2\subsetneq \ran V_1,$ we have
	\begin{equation}
		\overline{\cup_{z\in\D} \sigma(\varphi_1(z),\varphi_2(z))}=\sigma(M_{\varphi_1},M_{\varphi_2}).
	\end{equation}
	See Theorem \ref{pro:diag}, Theorem \ref{thm:negspec}, Theorem \ref{thm:pos:sp:mult} and Theorem \ref{offspecmulti}.
	
	In the following, when we consider the pair $(V_1,V_2)$ of commuting isometries, $V$  denotes the product $V_1V_2,$ and $\{e_n:n\in\Z\}$ denotes the standard orthonormal basis of $l^2(\Z).$
	
	We end this section with the comment that one of the strong points of the BCL theorem is that it models the $V_i$ in terms of functions of one variable whereas $V_i$ could, a priory, be dependent on two variables (multipliers on the Hardy space of the bidisc, for example). This strength will be greatly exploited in this note.

	\section{The defect operator}\label{sec:defopr}
	Recall that (\cite{Guo} and \cite{Yang}) the defect operator of a pair of commuting isometries $(V_1,V_2)$ is defined as 	
	\begin{equation}
		C(V_1,V_2)=I-V_1V_1^*-V_2V_2^*+V_1V_2V_2^*V_1^*.
    \end{equation}
		It is easy to see that (see \cite{Sarkar}) the defect
	\begin{equation}\label{eq:def:diffofproj}
		C(V_1,V_2) =P_{\ker V_1^*}-P_{V_2(\ker V_1^*)}=P_{\ker V_2^*}-P_{V_1(\ker V_2^*)},
	\end{equation}
	and \begin{equation}\label{eq:directsum}
		\ker V_1^*\oplus V_1(\ker V_2^*)=\ker V_2^*\oplus V_2(\ker V_1^*)=\ker V_1^*V_2^*.
	\end{equation}
	
	Note  that the  defect operator lives on $\ker V^*$ in the sense that the defect operator is zero on the orthogonal component of $\ker V^*.$  Equation \eqref{eq:def:diffofproj} shows that the defect is always a difference of two projections.  Let
	\begin{equation}\label{eq:P_1,P_2}
		P_1=P_{\ker V_1^*}\text{ and }P_2=P_{V_2(\ker V_1^*)}.
	\end{equation}
	Define
	\begin{align*}
		\H_1&:=\ran P_1\cap \ker P_2=\ker V_1^*\cap \ker V_2^*, \\ \H_2&:=\ran P_2\cap \ker P_1=V_1(\ker V_2^*)\cap V_2(\ker V_1^*),\\
		\H_3&:=\ran P_1\cap \ran P_2=\ker V_1^*\cap V_2(\ker V_1^*), \\
		\H_4&:=\ker P_1\cap \ker P_2=V_1(\ker V_2^*)\cap \ker V_2^*.
	\end{align*}
	Notice that $\H_i\perp \H_j$ if $i\ne j,$ and the $\H_i$ are reducing for $P_1$ and $P_2$ and  \begin{equation}\label{eq:directsumkerV^*}
		\H_1\oplus \H_2\oplus \H_3\oplus \H_4\subseteq \ker V^*.
	\end{equation}

	\begin{thm}\label{thm:iff}
		The following are equivalent:
	\begin{itemize}
		\item [(a)] 		The defect $C(V_1,V_2)$ is the difference of two mutually orthogonal projections.
		\item [(b)]	$\ker V^*=\H_1\oplus \H_2\oplus \H_3\oplus \H_4.$
		\item [(c)] $\ker V_1^*=\H_1\oplus \H_3.$
		\item [(d)] $V_1(\ker V_2^*)=\H_2\oplus \H_4.$
		\item [(e)] If $(\E,P,U)$ is the BCL triple for $(V_1,V_2),$
		\begin{equation}\label{eq:bcliff}
			U^*(\ran P)=(U^*(\ran P)\cap \ran P) \oplus (U^*(\ran P)\cap \ran P^\perp).
		\end{equation}
		
	\end{itemize}
	\end{thm}

\begin{proof}
Suppose
\begin{equation}\label{eq:Q1Q2}
	C(V_1,V_2)=Q_1-Q_2 \text{ with }  \ran Q_1\perp \ran Q_2
\end{equation}  for some projections $Q_1,Q_2$ in $\ker V^*.$ Notice that the pair $(Q_1,Q_2)$ satisfying \eqref{eq:Q1Q2} is {\em unique} if it exists. Then  for such a  pair $(V_1,V_2)$ we have
\begin{equation}\label{eq:diffprojn}
	C(V_1,V_2)=\begin{pmatrix}
		I_{\ran Q_1}&0&0\\0&-I_{\ran Q_2}&0\\0&0&0
	\end{pmatrix}=P_1-P_2.
\end{equation}
For any two projections $P$ and $Q$ in $\H$ and $x\in \H,$  $Px-Qx=x$ implies that $Px=x$ and $Qx=0.$ Using this fact and \eqref{eq:diffprojn}  we see that:
\begin{equation*}
	\H_1=\ran Q_1 \text{ and } \H_2=\ran Q_2.
\end{equation*}
Note that $P_1=P_2=\begin{pmatrix}
	I_\K&0\\0&0_\L
\end{pmatrix}$ on $\K\oplus \L$ for some $\K,\L$ that satisfy $\ker V^*= \H_1\oplus \H_2\oplus \K\oplus \L.$ Hence, by \eqref{eq:diffprojn} $\ran P_1=\H_1\oplus \K$ and $\ran P_2=\H_2\oplus \K.$ Thus  $\K=\ran P_1\cap \ran P_2=\H_3.$ Similarly $\L=\ker P_1\cap \ker P_2=\H_4.$ Therefore
\begin{equation*}
	\ker V^*=\H_1\oplus \H_2\oplus \H_3\oplus \H_4,
\end{equation*} and  in this decomposition
\begin{equation}\label{eq:matrixP_1,P_2}
	P_1=\begin{pmatrix}
		I&0&0&0\\0&0&0&0\\0&0&I&0\\0&0&0&0
	\end{pmatrix}, \quad
	P_2=\begin{pmatrix}
		0&0&0&0\\0&I&0&0\\0&0&I&0\\0&0&0&0
	\end{pmatrix}.
\end{equation}
This proves $(a)\implies (b),$ and  $(b)\implies (a)$ is trivial.

Note the fact that, if $\E_1,\E_2,\E_3$ and $\E_4$ are any subspaces of $\E$ satisfying $\E_1\oplus \E_2=\E=\E_3\oplus \E_4,$ then $\E_1=(\E_1\cap \E_3)\oplus(\E_1\cap \E_4)$ if and only if $\E_2=(\E_2\cap \E_3)\oplus(\E_2\cap \E_4).$ Using the above fact the equivalence $(c)\iff (d)$ follows from  \eqref{eq:directsum}.

The implication $(b)\implies (c)$ is clear. The implication $(c)\implies (b)$ follows from the equivalence of $(c)$ and $(d)$ and \eqref{eq:directsum}.

 Note that $\ker M_{\varphi_1}^*=\ran (I-M_{\varphi_1}M_{\varphi_1}^*)=1\otimes \ran (U^*PU)=1\otimes U^*(\ran P),$ $\ker M_{\varphi_2}^*=1\otimes \ran P^\perp$ and $M_{\varphi_2}(\ker M_{\varphi_1}^*)=  1\otimes \ran P.$ Now the proof of $(a)\iff (e)$ follows from the equivalence  $(a)\iff (c).$
\end{proof}
	Using the  $(e)$ part of the above theorem we easily get an example of a pair $(V_1,V_2)$ whose defect is not a difference of two mutually orthogonal projections.
\begin{example}\label{eg:outside}
	Let $\E=\C^2,U=\frac{1}{\sqrt{2}}\begin{pmatrix}
		1&1\\-1&1
	\end{pmatrix}$ and $P=\begin{pmatrix}
		1&0\\0&0
	\end{pmatrix}.$  The BCL triple $(\E,P,U)$ does not satisfy \eqref{eq:bcliff}, hence the defect of the pair $(V_1,V_2)$ corresponding to this BCL triple, is not the difference of two mutually orthogonal projections.
\end{example}

The Table \ref{classification:table} gives a neat classification. 
 We leave the proofs to the reader for the first two columns. The contents of the third column will unfold as we progress. One can notice that certain cases are not mentioned in the table. That is because those cases cannot occur.
	
%
%
%
%
%
	\begin{table}
		\begin{center}
			\caption{Classification}
	\begin{tabular}{|c |c|c|}
			\hline $\begin{array}{c}
					\ker V^*=\oplus_{i=1}^4\H_i \textbf{ and}\\				
			\textbf{the status of } \H_i		\end{array}$& {\bf The status of $(V_1,V_2)$}& {\bf BCL characterization}\\[1ex]
			\hline
			$\H_i=0 \text{ for }i=1,2,3,4$& $\begin{array}{c}
				C(V_1,V_2)=0 \text{ and}\\
				\text{ both } V_1\text{ and } V_2 \text{ are  unitaries}
			\end{array}$& $\E=0$\\
			\hline
			$\begin{array}{c}
				\H_i=0\text{ for } i=1,2,4,\\	\H_3\ne 0
			\end{array}	$
			& $\begin{array}{c}
				C(V_1,V_2)=0 \text{ and}\\
				V_2\text{ is a unitary}, V_1 \text{ is not a unitary}
			\end{array}$ & P=I \\
			\hline
			$\begin{array}{c}
				\H_i=0\text{ for } i=1,2,3,\\	\H_4\ne 0
			\end{array}	$
			&  $\begin{array}{c}
				C(V_1,V_2)=0 \text{ and}\\
				V_1\text{ is a unitary}, V_2 \text{ is not a unitary}
			\end{array}$ & P=0\\
			\hline
			$\begin{array}{c}
				\H_i=0\text{ for } i=1,2,\\	\H_i\ne 0 \text{ for } i=3,4
			\end{array}	$
			& $\begin{array}{c}
				C(V_1,V_2)=0 \text{ and}\\
				\text{ both } V_1\text{ and } V_2 \text{ are not unitaries}
			\end{array}$ & $\begin{array}{c}
		P \text{ is non-trivial and}\\
		U\text{ reduces } \ran P
		\end{array}$  \\
			\hline
			$\begin{array}{c}
				\H_1=0,\\	\H_i\ne 0 \text{ for } i=2,3,4
			\end{array}	$
			& $C(V_1,V_2)\le 0$ and $ C(V_1,V_2)\ne 0$ & $U(\ran P^\perp)\subsetneq \ran P^\perp$\\
			\hline
			$\begin{array}{c}
				\H_2=0,\\	\H_i\ne 0 \text{ for } i=1,3,4
			\end{array}$	
			& $C(V_1,V_2)\ge 0 $ and $C(V_1,V_2)\ne 0$ & $U(\ran P)\subsetneq \ran P$ \\
			\hline
			$\begin{array}{c}
				\H_i=0 \text{ for }i=3,4,\\	\H_i\ne 0 \text{ for } i=1,2
			\end{array}	$
			&   $\begin{array}{c}
				\ran V_1=\ran V_2  \text{ and}
				\\ V_i\text{ is not a unitary}
			\end{array}$ &$\begin{array}{c}
P\text{ is non-trivial and}\\
U(\ran P)=\ran P^\perp
		\end{array}$  \\	
			\hline 
			$\begin{array}{c}
				\H_4=0,\\	\H_i\ne 0 \text{ for } i=1,2,3
			\end{array}	$
			&  $\begin{array}{c}
				\ran V_1\subsetneq \ran V_2  \text{ and}
				\\ V_2 \text{ is not unitary}
			\end{array}$ & $\begin{array}{c}
		P\ne I \text{ and}	\\
        U(\ran P^\perp)\subsetneq\ran P
		\end{array}$ \\
			\hline
			$\begin{array}{c}
				\H_3=0,\\	\H_i\ne 0 \text{ for } i=1,2,4
			\end{array}	$
			&    $\begin{array}{c}
				\ran V_2\subsetneq \ran V_1 \text{ and}
				\\ V_1 \text{ is not unitary}
			\end{array}$ &$\begin{array}{c}
			P\ne 0 \text{ and}\\
			U(\ran P)\subsetneq\ran P^\perp			
		\end{array}$ \\
			\hline
			$\H_i\ne 0$ for $i=1,2,3,4$ & Unknown &--- \\[1ex] \hline
		\end{tabular}\label{classification:table}
			\end{center}	
	\end{table}

	\section{The Zero Defect Case} \label{ZeroDefectSection}
	
	\subsection{Structure}
	This subsection is mainly a rephrasing of known results.
	If one of the $V_i$'s is a unitary, then it is trivial to check that the defect $C(V_1, V_2)$ is zero. The following example is the prototypical example of a pure pair of commuting isometries with defect zero; this example serves as a building block in the general structure, see Theorem \ref{thm:zero struc}.
	
	\begin{example}\label{eg zero de}
		Let $\L$ be a non-zero Hilbert space and $W$ be a unitary on $\L.$ Consider the commuting pair of isometries $(M_{z}\otimes I_\L, I_{H^2_\D}\otimes W)$ on $H^2_{\D} \otimes \L$. As $I\otimes W$ is a unitary, the defect $C(M_{z}\otimes I_\L, I_{H^2_\D}\otimes W)=0.$ Also, $\sigma(M_{z}\otimes I_\L, I_{H^2_\D}\otimes W)=\Dc\times \sigma(W),$ by Lemma \ref{ReducingSpectrum} part \eqref{Part6ReducingSpectrum}. Also see \cite{Vasilescu}.
		
		Now consider the unitary $\Lambda:H^2_{\D}\otimes \L\to H^2_{\D}\otimes \L$ given by
		\begin{equation}
			\Lambda(\sum_{m= 0}^\infty a_mz^m)=\sum_{m= 0}^\infty W^m(a_m)z^m.
		\end{equation}
		Then, $\Lambda(M_z\otimes W^*)\Lambda^*=M_z\otimes I$ and $\Lambda (I\otimes W)\Lambda^*=I\otimes W.$ This says that $(M_z\otimes W^*,I\otimes W)$ and $(M_z\otimes I,I\otimes W)$ are jointly unitarily equivalent. In particular, we have $C(M_z\otimes W^*,I\otimes W)=0$ and
		\begin{equation}\label{eq:defzero:uni}
			\sigma (M_z\otimes W^*,I\otimes W)=\Dc\times \sigma(W)
		\end{equation}
		for any unitary $W.$
	\end{example}

	We proceed towards the structure of an arbitrary pair $(V_1, V_2)$ with $C(V_1, V_2) = 0$. 	A pair of commuting isometries $(V_1, V_2)$ is called {\em doubly commuting} if $V_1$ commutes with $V_2^*$. The following lemma is proved in \cite{Yang} and \cite{Sarkar}, it is relating positivity of the defect operator $C(V_1,V_2)$ with double commutativity of the pair $(V_1,V_2).$ We give a short proof here.    	
	\begin{lemma}\label{lem:SLOBCL}
		Let $(V_1,V_2)$ be a pair of commuting isometries on a Hilbert space $\H$. Then the following are equivalent:
		\begin{itemize}
			\item[(a)] $C(V_1,V_2)\geq 0.$
			\item[(b)] $(V_1,V_2)$ is doubly commuting.
			\item[(c)] If $(\E,P,U)$ is the BCL triple for $(V_1,V_2),$ then $U(\ran P)\subseteq \ran P.$
		\end{itemize}
	\end{lemma}
	
	\begin{proof}
		Since commuting unitaries are always doubly commuting, it is enough to prove the equivalences when $(V_1,V_2) = (M_{\varphi_1},M_{\varphi_2})$. By virtue of \eqref{DefectandBCL}, we have $C(V_1,V_2)\geq 0$ if and only if $U^*PU \ge P$ which happens if and only if $\ran P$ is invariant under $U.$ Now,
		\begin{align*}
			M_{\varphi_1}M_{\varphi_2}^*=M_{\varphi_2}^*M_{\varphi_1}&\text{ if and only if } (I-M_zM_z^*)\otimes (U^*PU^*P^\perp)=0 \\
			&\text{ if and only if } P^\perp UP=0\\
			&\text{ if and only if } \ran P \text{ is invariant under } U.
		\end{align*}
		This completes the proof.
	\end{proof}
	
	We recall some characterization results from \cite{Yang} and add a few new ones.	
	\begin{lemma}\label{lem:zeroequiv}
		Let $(V_1,V_2)$ be a pair of commuting isometries on a Hilbert space $\H.$ Then the following are equivalent:
		\begin{itemize}
			\item[(a)] $C(V_1,V_2)=0.$
			\item[(b)] $\ker V_2^*$ is a reducing subspace for $V_1$ and $V_1|_{\ker V_2^*}$ is a unitary.
			\item[(c)] $\ker V_1^*$ is a reducing subspace for $V_2$ and $V_2|_{\ker V_1^*}$ is a unitary.
			\item[(d)] The fringe operators $F_{1}$ and $F_{2}$ are unitaries.
			\item[(e)] $\ker V_1^*$ and $\ker V_2^*$ are orthogonal and their direct sum is $\ker V^*.$
			\item[(f)] $(\ran V_1\ominus \ran V)\oplus (\ran V_2\ominus \ran V)\oplus \ran V=\H.$
			\item[(g)] If $(\E,P,U)$ is the BCL triple for $(V_1,V_2),$ then $\ran P$ reduces $U.$

		\end{itemize}
		
		
	\end{lemma}
	
	\begin{proof}
		The equivalences of  $(a),$ $(b),$ $(c),$ $(d)$ and $(e)$ follows easily from   \eqref{eq:def:diffofproj} and \eqref{eq:directsum}.
		
		$(e)\Rightarrow (f)$: Suppose $(e)$ is true. We shall show that $\ker V_1^*=(\ran V_2\ominus \ran V).$ Suppose $x\in \ker V_1^*$, which implies $x\in \ran V_2$ and $x\in {(\ran V)}^{\perp}.$ So $x\in (\ran V_2\ominus \ran V).$
		If $x\in (\ran V_2\ominus \ran V),$ then $x\in \ran V_2$ and $x\in \ker V^*.$ So $x\in \ker V_1^*$. Similarly, $\ker V_2^*= \ran V_1\ominus \ran V.$ Hence $(f)$ is true.
		
		$(f) \Rightarrow (e)$: Suppose $(f)$ is true. Since $\ran V\subseteq \ran V_1,$ we have $\ran V_1= (\ran V_1\ominus \ran V)\oplus \ran V.$ Therefore $\ker V_1^*=(\ran V_1)^{\perp}=(\ran V_2\ominus \ran V).$ Similarly, $\ker V_2^*=(\ran V_1\ominus \ran V).$ Hence $\ker V_1^*\oplus \ker V_2^*=\ker V^*.$
		
		$(a) \Leftrightarrow (g)$: We use the formula  $C(V_1,V_2)=  (E_0\otimes (U^*PU-P))\oplus 0$ from \eqref{DefectandBCL}. This gives
		\begin{align*}
			C(V_1,V_2)=0 & \text{ if and only if } U^*PU - P = 0\\
			& \text{ if and only if } \ran P \text{ reduces } U.
		\end{align*}
		Thus, completes the proof.		
	\end{proof}
	
	We now write the structure theorem given in Popovici \cite[Sec. 4]{Popovici I}, which highlights the importance of Example \ref{eg zero de}. We give a proof for the completeness. 	
	\begin{thm} \label{thm:zero struc}
		Let $(V_{1},V_{2})$ be a pair of commuting isometries on  $\H$ with defect zero. Let $(\E,P,U)$ be the BCL triple for $(V_{1},V_{2}).$ Let $\E_1= \ran P$ and $\E_2 = \ran P^\perp $.   Then $\E_1,\E_2$ are reducing subspaces for $U$, i.e.,
		$$\E=\E_1\oplus\E_2, \ U=\begin{pmatrix}
			U_1 & 0 \\
			0 & U_2
		\end{pmatrix} \text{ and } P=\begin{pmatrix}
			I_{\E_1} & 0 \\
			0 & 0
		\end{pmatrix}\text{ in } \B(\E_1\oplus \E_2)$$	 for some unitaries $U_1$ and $U_2$ on $\E_1$ and $\E_2$ respectively.
		
		Also, $\H=(H^2_{\D}\otimes \E_1) \oplus (H^2_{\D}\otimes \E_2) \oplus \K$ where $\K=\underset{{n\ge 0}}{\bigcap} \ran (V_1V_2)^n$ and in this decomposition, $$V_1= \begin{pmatrix}
			M_z\otimes I_{\E_1}& 0&0\\
			0& I_{H^2_{\D}}\otimes U_2^* & 0\\
			0&0&W_1
		\end{pmatrix}, \quad V_2= \begin{pmatrix}
			I_{H^2_{\D}}\otimes U_1 & 0&0\\
			0& M_z\otimes I_{\E_2} & 0\\
			0&0&W_2
		\end{pmatrix},$$  up to unitarily equivalence, for some unitary $U_i$ on $\E_i, i=1,2$ and commuting unitaries $W_1, W_2$ on $\K.$
	\end{thm}
	
	\begin{proof}
		Let $(\E,P,U)$ be the BCL triple for $(V_1,V_2).$ By Lemma \ref{lem:zeroequiv}, we have $\E_1$ reduces $U$ and in the decomposition $\E=\E_1\oplus\E_2,$ we have $$ U=\begin{pmatrix}
			U_1 & 0 \\
			0 & U_2
		\end{pmatrix} \text{ and } P=\begin{pmatrix}
			I_{\E_1} & 0 \\
			0 & 0
		\end{pmatrix}.$$
		In this case,                        $\varphi_1(z)=\begin{pmatrix}
			zU_1^{*} & 0 \\
			0 & U_2^{*}
		\end{pmatrix}$ and $\varphi_2(z)=\begin{pmatrix}
			U_1 & 0 \\
			0 & zU_2
		\end{pmatrix},$ for $z\in \D$. Therefore, $$M_{\varphi_1}=\begin{pmatrix}M_z\otimes U_1^*& 0\\
			0& I_{H^2_{\D}}\otimes U_2^*\end{pmatrix},\  M_{\varphi_2}=\begin{pmatrix}I_{H^2_{\D}}\otimes U_1& 0\\
			0& M_z\otimes U_2\end{pmatrix}.$$
		Note that $(M_z\otimes U_1^*,I_{H^2_{\D}}\otimes U_1)$  and $(M_z\otimes I_{\E_1},I_{H^2_{\D}}\otimes U_1)$ are jointly unitarily equivalent and $(I_{H^2_{\D}}\otimes U_2^*,M_z\otimes U_2)$ and $(I_{H^2_{\D}}\otimes U_2^*,M_z\otimes I_{\E_2})$ are jointly unitarily equivalent; see Example \ref{eg zero de}. This completes the proof, by Theorem \ref{thm:BCL}.
	\end{proof}
	
	\begin{remark}\label{rem:intermsofVi}
		In the structure theorem (Theorem \ref{thm:zero struc}) one can write the $U_1, U_2$ and $\E_1,\E_2$ explicitly in terms of $V_1$ and $V_2$ as follows: $\E_1=\ker V_1^*, \E_2=\ker V_2^*$ and	
		$U_1=V_2|_{\ker V_1^*},U_2=V_1^*|_{\ker V_2^*}.$ 
		This is because of Lemma \ref{lem:zeroequiv}.
	\end{remark}

	A comment is in order. Over several decades Marek S{\l}oci\'{n}ski, first by himself and then with his collaborators, has developed a complete structure theorem on commuting pairs of isometries; see \cite{BKPS} and the references therein. One of his early results is the following.

	\begin{thm}[M. S{\l}oci\'{n}ski \cite{Slocinski}]\label{Slocinski}
		Let  $(V_1, V_2)$ be a pair of doubly commuting isometries on
		a Hilbert space $\H.$ Then there exists a unique decomposition
		$$\H = \H_{ss}\oplus \H_{su}\oplus \H_{us}\oplus \H_{uu},$$
	where the subspace $\H_{ij}$ reduces both $V_1$ and $V_2$ for all $i, j \in \{s, u\}$. Moreover, $V_1$ on $\H_{ij}$ is a
		shift if $i = s$ and unitary if $i = u$ and $V_2$ is a shift if $j = s$ and unitary if $j = u$.
	\end{thm}
By Theorem \ref{thm:zero struc}, and the fact that if one of the $V_i$'s is a unitary then  the defect is zero, we have $C(V_1,V_2)=0$ if and only if $(V_1,V_2)$ is doubly commuting and $\H_{ss}$ in Theorem \ref{Slocinski} is $\{0\}.$
	
	\subsection{Joint spectrum}
	If $(V_1,V_2)$ is pure and defect $C(V_1,V_2)=0,$ then by Theorem \ref{thm:zero struc} and Remark \ref{rem:intermsofVi},
	\begin{equation}
		\sigma(V_1,V_2)=\begin{cases}
			\Dc\times \sigma(U_1) &\text{if } V_2 \text{ is unitary,}\\
			\sigma(U_2^*)\times\Dc &\text{if } V_1 \text{ is unitary,}\\
			\Dc\times \sigma (U_1)\cup \sigma(U_2^*)\times \Dc & \text{if neither }V_1 \text{ nor } V_2 \text{ is a unitary,}
		\end{cases}
	\end{equation}
	where $U_1=V_2|_{\ker V_1^*},U_2=V_1^*|_{\ker V_2^*}.$
	
	
	\begin{lemma}\label{lem:zero}
		Let $(V_1,V_2)$ be a pair of commuting isometries with defect zero and $\ker V^*\ne \{0\}.$ Let $(\E,P,U)$ be the BCL triple for $(V_1,V_2).$  Let $U_1=U|_{\ran P}$ and $U_2=U|_{\ran P^\perp}.$ Then
		\begin{equation}
			\sigma(\varphi_1(z),\varphi_2(z))=
			\begin{cases}
				\{(z\overline{\lambda},\lambda): \lambda\in \sigma(U_1)\}& \text{ if } V_2 \text{ is  a unitary},\\
				\{(\overline{\mu},z\mu):\mu\in \sigma(U_2)\}& \text{ if } V_1 \text{ is  a unitary,}\\
				\{(z\overline{\lambda},\lambda), (\overline{\mu},z\mu): \lambda\in \sigma(U_1),\ \mu\in \sigma(U_2)\}&   \text{ otherwise.}
			\end{cases}
		\end{equation}
		Here, for every point in the joint spectrum, the non-singularity breaks at stage 3.	
	\end{lemma}

	\begin{proof}
		The proof in the case when neither $V_1$ nor $V_2$ is a unitary is done below in detail. The other cases follow similarly.
		
		Letting $\E_1=\ran P$ and $\E_2=\ran P^\perp$, it is an easy check that both $\E_1$ and $\E_2$ are non-trivial. By Lemma \ref{lem:zeroequiv}, $\E_1$ reduces $U.$ Hence in the decomposition $\E=\E_1\oplus \E_2,$ we have
		$$\varphi_1(z)=\begin{pmatrix}
			zU_1^{*} & 0 \\
			0 & U_2^{*}
		\end{pmatrix} \text{ and } \varphi_2(z)=\begin{pmatrix}
			U_1 & 0 \\
			0 & zU_2
		\end{pmatrix},$$
		for $z\in \D.$		Then,
		$$\sigma(\varphi_1(z))=\sigma(zU_1^{*}) \cup \sigma(U_2^{*}) \text{ and } \sigma(\varphi_2(z))=\sigma(U_1) \cup \sigma (z U_2) $$ for all $z\in \D$. For $z\in \D,$ we have
		$$  \varphi_1(z)\varphi_2(z)=zI_{\E_1\oplus \E_2}=\varphi_2(z)\varphi_1(z).$$
		Therefore, for all $z\ne 0,$ by Lemma \ref{ReducingSpectrum} part \eqref{Part1ReducingSpectrum} and polynomial spectral mapping theorem, we get
		$$\sigma(\varphi_1(z),\varphi_2(z))\subseteq \{(z\overline{\lambda},\lambda), (\overline{\mu},z\mu): \lambda\in \sigma(U_1),\ \mu\in \sigma(U_2)\}.$$ For $z=0,$ it is easy to see that$$\sigma(\varphi_1(0),\varphi_2(0))\subseteq \{(0,\lambda), (\overline{\mu},0): \lambda\in \sigma(U_1),\ \mu\in \sigma(U_2)\}.$$
		
		To prove the other containments, let $\lambda\in \sigma(U_1).$ Then $U_1-\lambda I$ is not onto, because $U_1$ is normal. Let
		$h=\begin{pmatrix}
			h_1\\
			h_2
		\end{pmatrix}, k=\begin{pmatrix}
			k_1\\
			k_2
		\end{pmatrix}\in \E_1\oplus \E_2. $
		Then,
		$$(\varphi_1(z)-z\overline{\lambda}I)k+(\varphi_2(z)-{\lambda}I)h=\begin{pmatrix}
			z(U_1^{*}-\overline{\lambda}I) k_1+ (U_1-\lambda I)h_1 \\
			(U_2^{*}-z \overline{\lambda}I)k_2+(zU_2-\lambda I)h_2
		\end{pmatrix}.$$
		Since $U_1$ is a normal operator,
		$\ran(U_1-\lambda I)=\ran(U_1^{*}-\overline{\lambda} I)$
		and hence the first component of the above spans only $\ran(U_1-\lambda I)$. Hence we have
		$$\ran(\varphi_1(z)-z\overline{\lambda}I)+\ \ran(\varphi_2(z)-{\lambda}I)\neq \E_1\oplus \E_2.$$ Therefore $(z\overline{\lambda},\lambda)\in \sigma(\varphi_1(z),\varphi_2(z))$.
		
		Similarly, if $\mu\in \sigma(U_2),$ we have
		$$\ran(\varphi_1(z)-\overline{\mu}I)+\ \ran(\varphi_2(z)-z{\mu}I)\neq \E_1\oplus \E_2.$$ Therefore $(\overline{\mu},z\mu)\in \sigma(\varphi_1(z),\varphi_2(z)).$ So,
		\begin{equation*}
			\sigma(\varphi_1(z),\varphi_2(z))=\{(z\overline{\lambda},\lambda), (\overline{\mu},z\mu): \lambda\in \sigma(U_1),\ \mu\in \sigma(U_2)\} \text{ for } z\in \D.
		\end{equation*}
	\end{proof}

	\begin{thm}\label{pro:diag}
		
		
		With the hypothesis of Lemma \ref{lem:zero}, we also have
		\begin{equation*}
			\sigma(M_{\varphi_1},M_{\varphi_2})=\overline{\cup_{z\in \D}\sigma(\varphi_1(z),\varphi_2(z))}=
			\begin{cases}
				\Dc\times \sigma(U_1) & \text{ if } V_2 \text{ is  a unitary,}\\
				\sigma(U_2^{*})\times \Dc & \text{ if } V_1 \text{ is  a  unitary,}\\
				\Dc\times \sigma(U_1)\cup \sigma(U_2^{*})\times \Dc & \text{ otherwise.}
			\end{cases}
		\end{equation*}
		and, for every point $(z_1,z_2)\in \begin{cases}
			\D\times \sigma(U_1) & \text{ if }  V_2 \text{ is  a unitary,}\\
			\sigma(U_2^{*})\times \D & \text{ if } V_1 \text  { is  a unitary,}\\
			\D\times \sigma(U_1)\cup \sigma(U_2^{*})\times \D & \text{ otherwise,}
		\end{cases}$\\
		the non-singularity of $K(M_{\varphi_1}-z_1I,M_{\varphi_2}-z_2I)$ breaks at stage 3.
	\end{thm}
	
	\begin{proof}
		%
		
		As in Lemma \ref{lem:zero}, we prove only the case when neither $V_1$ nor $V_2$ is a unitary, other cases follow similarly.
		
		We saw in the proof of Lemma \ref{lem:zero} that for any point $z\in \D$, a pair of points $(z_1,z_2)\in \sigma(\varphi_1(z),\varphi_2(z))$ if and only if
		$$\ran(\varphi_1(z)-z_1I)+\ \ran(\varphi_2(z)-z_2I)\neq \E_1\oplus \E_2.$$
		Hence, if $(z_1,z_2)\in \sigma(\varphi_1(z),\varphi_2(z)),$
		\begin{equation}\label{eq:chainbreak}
			\ran(M_{\varphi_1}-z_1I)+\ \ran(M_{\varphi_2}-z_2I)\neq H^{2}(\E_1\oplus \E_2).
		\end{equation}
		Therefore, $(z_1,z_2)\in \sigma(M_{\varphi_1},M_{\varphi_2}),$  which implies that
		\begin{equation*}
			{\cup_{z\in \D}\sigma(\varphi_1(z),\varphi_2(z))}\subseteq\sigma(M_{\varphi_1},M_{\varphi_2}).
		\end{equation*}
		Note that
		\begin{align}\label{eq:inc3}
			{\cup_{z\in \D}\sigma(\varphi_1(z),\varphi_2(z))} & =\cup_{z\in \D}\{(z\overline{\lambda},\lambda), (\overline{\mu},z\mu): \lambda\in \sigma(U_1),\ \mu\in \sigma(U_2)\}\nonumber\\
			& ={\D}\times \sigma(U_1)\cup \sigma(U_2^{*})\times {\D}.\end{align}
		Since
		$$M_{\varphi_1}=\begin{pmatrix}
			M_{z}\otimes U_1^{*} & 0 \\
			0 & I_{H^2_\D}\otimes U_2^{*}
		\end{pmatrix} \text{ and } M_{\varphi_2}=\begin{pmatrix}
			I_{H^2_\D}\otimes U_1 & 0 \\
			0 & M_{z}\otimes U_2
		\end{pmatrix},$$ we have
		$$\sigma(M_{\varphi_1},M_{\varphi_2})=\sigma(M_{z}\otimes U_1^{*}, I_{H^2_\D}\otimes U_1) \cup \sigma(I_{H^2_\D}\otimes U_2^{*}, M_{z}\otimes U_2).$$
		By \eqref{eq:defzero:uni}, we have
		\begin{equation}\label{eq:inc2}
			\sigma(M_{\varphi_1},M_{\varphi_2})= \Dc\times \sigma(U_1)\cup \sigma(U_2^{*})\times \Dc.
		\end{equation}
		Therefore from \eqref{eq:inc3} and \eqref{eq:inc2} we have
		\begin{equation*}
			\sigma(M_{\varphi_1},M_{\varphi_2})=\Dc\times \sigma(U_1)\cup \sigma(U_2^{*})\times \Dc=\overline{\cup_{z\in \D}\sigma(\varphi_1(z),\varphi_2(z))}.
		\end{equation*}
		The final thing to note is that for every point $(z_1,z_2)$ in the set $\D\times \sigma(U_1)\cup \sigma(U_2^{*})\times {\D},$ the non-singularity breaks at stage 3 and this is a direct consequence of \eqref{eq:chainbreak}.
	\end{proof}
	
	To conclude the section, we note that by Theorem \ref{thm:BCL}, $\sigma(V_1,V_2)=\sigma(M_{\varphi_1},M_{\varphi_2})\cup \sigma(V_1|_{\H_u},V_2|_{\H_u}).$ Hence Theorem \ref{pro:diag} tells that
	\begin{equation}\label{eq:union}
		\overline{\cup_{z\in \D}\sigma(\varphi_1(z),\varphi_2(z))}\subseteq \sigma(V_1,V_2).
	\end{equation}
	The equality in \eqref{eq:union} holds if and only if $\sigma(V_1|_{\H_u},V_2|_{\H_u})\subseteq \sigma(M_{\varphi_1},M_{\varphi_2}).$

	\section{The Negative Defect Case}

	\subsection{The Prototypical Example} \label{subsec:neg:eq}

	When the defect operator is negative, a fundamental example plays an important role in much the same way the unilateral shift plays its role in the Wold decomposition of a single isometry. We shall first describe this example and then show how it is a part of every pair of commuting isometries with negative defect.
	
	The {\em Hardy space of $\mathcal E$-valued functions on the bidisc $\mathbb D^2$} is
	\begin{multline*} H^2_{\mathbb D^2}({\mathcal E})  =\{f : \mathbb D^2 \rightarrow \mathcal E \mid f \text{ is analytic and } f(z_1,z_2) = \sum_{m,n=0}^{\infty} a_{m,n} z_1^mz_2^n \\\text{ with } \sum_{m,n=0}^{\infty} \|a_{m,n}\|_\E^2 < \infty\}.\end{multline*}
	This is a Hilbert space with the inner product
	$$\langle \sum_{m,n=0}^{\infty} a_{m,n} z_1^mz_2^n , \sum_{m,n=0}^{\infty} b_{m,n} z_1^mz_2^n \rangle =  \sum_{m,n=0}^{\infty} \langle a_{m,n} , b_{m,n} \rangle_\E$$ and is identifiable with $H^2_{\mathbb D^2} \otimes \mathcal E$ where $H^2_{\mathbb D^2}$ stands for the Hardy space of scalar-valued functions on $\mathbb D^2$.
	
	Let $U:H^2_{\D^2}\to H^2_{\D^2}$ be the unitary defined by
	\begin{equation} \label{specialunitary}
		U(z_1^{m_1}z_2^{m_2})=
		\begin{cases}
			z_1^{m_1+2}z_2^{m_2} & \text{if } m_1\ge m_2,\\
			z_1^{m_1+1}z_2^{m_2-1} & \text{if }  m_1+1=m_2,\\
			z_1^{m_1}z_2^{m_2-2} & \text{if }  m_1+2\le m_2.
		\end{cases}
	\end{equation}
	on the orthonormal basis $\{z_1^{m_1}z_2^{m_2}\}_{m_1,m_2\geq  0}.$
	
	The pair of multipliers by the coordinate functions $(M_{z_1}, M_{z_2})$ forms a pair of doubly commuting isometries on $H^2_{\mathbb D^2}$. There is a natural isomorphism between the Hilbert spaces $H^2_{\mathbb D^2}$ and $H^2_{\mathbb D} \otimes H^2_{\mathbb D}$ wherein $z_1^{m_1} z_2^{m_2}$ is identified with $z_1^{m_1} \otimes z_2^{m_2}$. In this identification, the pair of coordinate multipliers $(M_{z_1}, M_{z_2})$ is identified with $(M_z \otimes I, I \otimes M_z)$.
	
	\begin{definition}
		
		The pair of bounded operators $\tau_1:=U^*M_{z_1}$ and $\tau_2:=M_{z_2}U$ on the Hardy space of the bidisc $H^2_{\mathbb D^2}$ will be called the fundamental isometric pair with negative defect.
	\end{definition}
	
	The following lemma justifies the name except the word $fundamental$ which will be clear from Theorem \ref{StructureNegative}.
	
	\begin{lemma}
		The pair $(\tau_1,\tau_2)$ is a pair of commuting isometries with  defect negative and non-zero.
	\end{lemma}
	
	\begin{proof}
		It is simple to check that the unitary $U$ defined in \eqref{specialunitary} commutes with $M_{z_1z_2}$, the operator of multiplication by the function $z_1z_2$. That proves commutativity of $\tau_1$ and $\tau_2$. They are isometries because each is a product of an isometry and a unitary. Now, let
		$$\W_1=\ker(\tau_1^*)=\spncl\{z_2^2,z_2^3,z_2^4,\dots \} \text{ and } \W_2=\ker(\tau_2^*)=\spncl\{1,z_1,z_1^2,\dots \}. $$
		Then, $\tau_2(\W_1)=\spncl\{z_2,z_2^2,z_2^3,\dots\}.$ Thus, we have
		\[C(\tau_1,\tau_2)=P_{\W_1}-P_{\tau_2(\W_1)}=-P_{\spn \{z_2\}}\le 0 .\]
		That completes the proof.
	\end{proof}
	
	We shall prove the following lemma to compute the joint spectrum of $(\tau_1,\tau_2).$
	
	\begin{lemma} \label{lem:neg}
		For any $\lambda\in \D$, we have
		\begin{equation}\label{eq:lem}
			(\ran(\tau_2-\lambda I))^\perp
			=\left\{ (I - \overline{\lambda} \tau_2)^{-1} x : x \in \ker M_{z_2}^*\right\}.
		\end{equation}
	\end{lemma}

\begin{proof}
	Using the Neumann series $(I-\bar{\lambda} A)^{-1}=\sum_{n\ge 0}\bar{\lambda}^nA^n$ for  $\lambda \in \D$ and any contraction $A,$ it is straightforward that  the equality $$(\ran (A-\lambda I))^\perp=\{(I-\bar{\lambda}A)^{-1}x:x\in \ker A^*\}$$ is satisfied when $A$ is an isometry. The proof is complete by noting that $\ker \tau_2^*=\ker U^*M_{z_2}^*=\ker M_{z_2}^*.$
\end{proof}
	
%
%
	
	Recall the Koszul complex for a pair of commuting bounded operators $(T_1, T_2)$ from \eqref{eq:Kozul}:
	\begin{equation}\label{eq:KozulRepeat}
		0\overset{\delta_0}{\to} \H \overset{\delta_1}{\to}\H\oplus \H\overset{\delta_2}{\to}\H\overset{\delta_3}{\to} 0.
	\end{equation}
	It is well-known that the most difficult stage to treat for the purpose of showing lack of exactness is the stage 2.
	
	\begin{proposition}\label{pro:eg}
		The fundamental isometric pair with negative defect has the full closed bidisc $\overline{\D^2}$ as its joint spectrum. Moreover, for every point in the open bidisc $\D^2,$ the non-singularity breaks at stage 2.
	\end{proposition}
	
	\begin{proof}
		Let $\lambda_1,\lambda_2\in \D.$ We shall find a non-zero function $ h_2\in (\ran(\tau_2-\lambda_2I))^\perp$ such that
		\begin{equation} \label{taulambdah2}(\tau_1-\lambda_1I)h_2\in \ran(\tau_2-\lambda_2I).\end{equation}
		This would imply that there exists $h_1\in \H=H^2_{\D ^2}$ such that
		\[(\tau_1-\lambda_1I)h_2=(\tau_2-\lambda_2I)h_1\]
		producing a pair $(h_1,h_2)$ in $\ker \delta_2$ which would not be in $\ran\delta_1$.
		
		To that end,  we shall use the description of $(\ran(\tau_2-\lambda I))^\perp$ obtained in Lemma \ref{lem:neg}.
		Since any element from $\ker \tau_2^*=\ker M_{z_2}^*$ is of the form $\sum_{m=0}^\infty a_mz_1^m $ for a square summable sequence $\{a_m\}_{m\ge 0}$, it follows from Lemma \ref{lem:neg} that
		\begin{align}
			(\ran(\tau_2-\lambda_2I))^\perp
			&=\{ (I - \overline{\lambda_2} \tau_2)^{-1} \sum_{m= 0}^\infty a_mz_1^m : \sum_{m= 0}^\infty |a_m|^2 < \infty\} \nonumber \\
			&=\{\sum_{n= 0}^\infty(\overline{\lambda_2}M_{z_2} U)^n (\sum_{m= 0}^\infty a_mz_1^m):\sum_{m= 0}^\infty |a_m|^2 < \infty\} \nonumber \\
			&=\{\sum_{m,n=0}^{\infty}\overline{\lambda_2}^na_mz_1^{m+2n}z_2^n:\sum_{m= 0}^\infty |a_m|^2 < \infty\}. \label{StrangeSeries}
		\end{align}		
		
		Our candidate for $h_2$ to satisfy \eqref{taulambdah2} is
		$$h_2= \frac{1}{(1 - \overline{\lambda_2}z_1^2z_2)(1 - \lambda_1z_1)} = \sum_{m,n=0}^{\infty}\overline{\lambda_2}^n \lambda_1^m z_1^{m+2n}z_2^n.$$
		By \eqref{StrangeSeries}, this function is in $(\ran(\tau_2-\lambda_2I))^\perp$.
		We shall verify below that $(\tau_1-\lambda_1I)h_2$ is in the closure of $\ran(\tau_2-\lambda_2I)$. Since $|\lambda_2|<1$ and $\tau_2$ is an isometry, $(\tau_2-\lambda_2 I)$ is bounded below and hence its range is closed. That will complete the proof.
		
		First note that
		\begin{equation}\label{ActionofU}
			U \sum_{m,n=0}^{\infty} {\overline{\lambda_2}}^n a_m z_1^{m+2n} z_2^n
			=  \sum_{m,n=0}^{\infty} {\overline{\lambda_2}}^n a_m z_1^{m+2n+2} z_2^n
			=  M_{z_1}^2 \sum_{m,n=0}^{\infty} {\overline{\lambda_2}}^n a_m z_1^{m+2n} z_2^n
		\end{equation}
		and
		\begin{align}
			(M_{z_1}^*-\lambda_1 I)h_2
			= & M_{z_1}^*(h_2)-\lambda_1 h_2 \nonumber \\
			= & \underset{(m,n)\ne (0,0)}{\sum_{m,n=0}^{\infty}}\overline{\lambda_2}^n \lambda_1^m z_1^{m+2n-1}z_2^n-\sum_{m,n=0}^{\infty} \overline{\lambda_2}^n \lambda_1^{m+1} z_1^{m+2n}z_2^n \nonumber \\
			= & \sum_{m=1,n=0}^{\infty}\overline{\lambda_2}^n \lambda_1^m z_1^{m+2n-1}z_2^n+ \sum_{n=1}^{\infty}\overline{\lambda_2}^n z_1^{2n-1}z_2^n-\sum_{m,n=0}^{\infty} \overline{\lambda_2}^n \lambda_1^{m+1} z_1^{m+2n}z_2^n \nonumber \\
			= & \sum_{m,n=0}^{\infty} \overline{\lambda_2}^n \lambda_1^{m+1} z_1^{m+2n}z_2^n + \sum_{n=1}^{\infty}\overline{\lambda_2}^n z_1^{2n-1}z_2^n-\sum_{m,n=0}^{\infty} \overline{\lambda_2}^n \lambda_1^{m+1} z_1^{m+2n}z_2^n \nonumber \\
			= & \sum_{n=1}^{\infty}\overline{\lambda_2}^n z_1^{2n-1}z_2^n .
		\end{align}
		We now compute the inner product between a typical element of $(\ran(\tau_2-\lambda_2I))^\perp$ and $(\tau_1-\lambda_1I)h_2$ by using the two equations above.
		\begin{align*}
			&\la (\tau_1-\lambda_1I)h_2, \sum_{m,n=0}^{\infty} \overline{\lambda_2}^{n}a_{m}z_1^{m+2n}z_2^{n}\ra \\
			&=\la (U^*M_{z_1}-\lambda_1I)h_2, \sum_{m,n=0}^{\infty} \overline{\lambda_2}^{n}a_{m}z_1^{m+2n}z_2^{n} \ra \\
			&=\la M_{z_1}h_2, U\sum_{m,n=0}^{\infty} \overline{\lambda_2}^{n}a_{m}z_1^{m+2n}z_2^{n} \ra- \lambda_1\la h_2,\sum_{m,n=0}^{\infty} \overline{\lambda_2}^{n}a_{m}z_1^{m+2n}z_2^{n} \ra\\
			&= \la M_{z_1}h_2,  M_{z_1}^2\sum_{m,n=0}^{\infty} \overline{\lambda_2}^{n}a_{m}z_1^{m+2n}z_2^{n} \ra- \lambda_1\la h_2,\sum_{m,n=0}^{\infty} \overline{\lambda_2}^{n}a_{m}z_1^{m+2n}z_2^{n} \ra\\
			&= \la (M_{z_1}^*-\lambda_1 I)h_2, \sum_{m,n=0}^{\infty} \overline{\lambda_2}^{n}a_{m}z_1^{m+2n}z_2^{n} \ra\\
			&= \la	\sum_{n=1}^{\infty}\overline{\lambda_2}^n z_1^{2n-1}z_2^n, \sum_{m,n=0}^{\infty} \overline{\lambda_2}^{n}a_{m}z_1^{m+2n}z_2^{n} \ra
			= 0.
		\end{align*}
		This shows that $(\tau_1-\lambda_1I)h_2\in \overline{\ran}(\tau_2-\lambda_2I) = \ran(\tau_2-\lambda_2I)$  and hence  completes the proof.
	\end{proof}

	\begin{note}In Remark \ref{NrestrictedP}, we shall see that there is a joint invariant subspace $\mathcal M$ for $(\tau_1, \tau_2)$ such that the defect operator of $(\tau_1|_{\mathcal M}, \tau_2|_{\mathcal M})$ is positive (and not zero). \end{note}

	\subsection{General Theory for the negative defect case:} Here we shall show that the fundamental example above is a typical example. This helps us to compute the joint spectrum of any commuting pair of isometries with negative defect.
	
		In \cite{GG}, Gaspar and Gaspar  introduced the dual doubly commuting pairs. If $(\bar{V_1},\bar{V_2})$ is the {\em minimal unitary extension} to $\bar{\H}$ of $(V_1,V_2)$ acting on $\H$, then the pair of commuting isometries $(\bar{V_1}^*|_{\bar \H \ominus \H},  \bar{V_2}^*|_{\bar \H \ominus \H})$ is called the {\em dual} of $(V_1,V_2)$. If the dual is doubly commuting, then  $(V_1,V_2)$ is called a {\em dual doubly commuting} pair.
		
		A pair $(V_1,V_2)$ of commuting isometries is called a {\em bi-shift} (see \cite{Popovici II}) if there is a wandering subspace $\mathcal R$ (i.e., $ V_1^{p_1}V_2^{p_2}(\mathcal R)\perp V_1^{q_1}V_2^{q_2}(\mathcal R)$ if $(p_1, p_2), (q_1, q_2)\in \Z_+^2$ and $(p_1, p_2) \ne (q_1, q_2)$) such that
		$$\H=\bigoplus_{(n_1, n_2) \in \Z_+^2 } V_1^{n_1}V_2^{n_2}(\mathcal R).$$
		   $(V_1,V_2)$ is called a {\em modified bi-shift} if it is pure and its dual is a bi-shift.	
		
		 Popovici used the concepts above greatly in his papers \cite{Popovici II} and \cite{Popovici}. We are thankful to him for sending us his papers. First we shall give a characterizing lemma for this case;  see also \cite{Yang}.
	\begin{lemma}\label{lem:UP}
		Let $(V_1,V_2)$ be a pair of commuting isometries on a Hilbert space $\H$. Then the following are equivalent:
		\begin{itemize}
			\item[(a)] $C(V_1,V_2)\leq 0$ and $C(V_1,V_2)\neq 0.$
			\item[(b)] $V_2(\ker V_1^*)\supsetneq \ker V_1^*.$
			\item[(c)] $V_1(\ker V_2^*)\supsetneq \ker V_2^*.$
			\item[(d)] The adjoint of the fringe operators are isometries and not unitaries.
				\item[(e)] $(V_1,V_2)$ is dual doubly commuting and $C(V_1,V_2)\neq 0.$			
			\item[(f)] $\ker V_1^*$ is orthogonal to $\ker V_2^*$ and $\ker V_1^*\oplus \ker V_2^*\neq \ker V^*.$
			\item[(g)] $C(V_1,V_2)$ is the negative of a non-zero projection.
			\item[(h)] If $(\E,P,U)$ is the BCL triple for $(V_1,V_2),$ then $U(\ran P^\perp)\subsetneq \ran P^\perp.$
		\end{itemize}
	\end{lemma}
	\begin{proof}
		The equivalence $(e) \Leftrightarrow (h)$ is proved in \cite{GG}. All other proof are along the same lines as the proofs of various parts of Lemma \ref{lem:zeroequiv}.	
	\end{proof}

The geometrical structure and a model for dual doubly commuting isometries is known due to \cite{GG,Popovici I}. Here we observe that the multiplication operators $M_{\varphi_i},i=1,2$ associated to $(V_1,V_2)$ as in Theorem \ref{thm:BCL}, have some  special forms in this case. This also helps us getting a model in the Hardy space of the bidisc. Some steps of the proof are used to obtain Theorem \ref{thm:negspec}. The wandering space arguments used in the proof of the following theorem,  appears in \cite[Thm 4.3]{Popovici I}.

 Let $\psi_1,\psi_2:\D\to \B(l^2(\Z))$ be the multipliers associated with the BCL triple $(l^2(\Z),p_-,\omega)$ where $p_-$ is  the projection onto $\spncl\{e_n:n<0\}$ and $\omega$ is the bilateral shift on $l^2(\Z)$.

	\begin{thm}\label{thm:NegStructure}
		Let $(V_1,V_2)$ be a pair of commuting isometries such that $C(V_1,V_2)$ is a non-zero negative operator. Let  $(\E,P,U)$ be the  BCL triple for $(V_1,V_2).$ Then, up to unitary equivalence $$\E = (l^2(\Z)\otimes \L) \oplus \E_2$$
		for some non-trivial closed subspace $\L$ and a closed subspace $\E_2$ of $\E.$ Moreover,
		\begin{equation}\label{eq:mphimpsi} M_{\varphi_i}=\begin{pNiceMatrix}[first-row,last-col,nullify-dots]
				H^2_\D(l^2(\Z)) \otimes \L& H^2_\D(\E_2)&\\
				M_{\psi_i}\otimes I_\L & 0&H^2_\D(l^2(\Z)) \otimes \L \\ 0& {M_{\varphi_i}}|_{H^2_{\D}(\E_2)}& H^2_\D(\E_2)
			\end{pNiceMatrix}
		\end{equation}
		with $C({M_{\varphi_1}}|_{H^2_{\D}(\E_2)},{M_{\varphi_2}}|_{H^2_{\D}(\E_2)})=0.$ In particular,
		\begin{equation}\label{eq:incpsi}
			\sigma(M_{\psi_1},M_{\psi_2})= \sigma( M_{\psi_1}\otimes I_\L,M_{\psi_2}\otimes I_\L)\subseteq\sigma(M_{\varphi_1},M_{\varphi_2})\subseteq \sigma(V_1,V_2).
		\end{equation}
	\end{thm}
	
	\begin{proof}
		Since $C(V_1,V_2)\le 0 \text{ and }C(V_1,V_2)\ne 0$,  by Lemma \ref{lem:UP}, $U(\ran P^\perp ) \subsetneq  \ran P^\perp.$ Consider $\L:=(\ran P^\perp \ominus U(\ran P^\perp))\neq \{0\}.$ Let $m,n\in \Z$ and $m> n$. Then for $x,y\in \L$ we have \[\la U^m x, U^ny\ra=\la U^{m-n}x, y\ra=0,\] because $U^{m-n}x\in U(\ran P^\perp).$ 
		Therefore $U^m(\L)\perp U^n(\L)$ if $ m, n\in \Z$ and $m\ne n.$
		
		Set $\E_1:=\oplus_{n\in \Z} U^n(\L).$ Clearly $U$ reduces $\E_1.$ Now  $U(\ran P^\perp )\subsetneq \ran P^\perp$ implies:
		\begin{equation}\label{eq:n>0}
			\oplus_{n\ge 0}U^n(\L)\subseteq \ran P^\perp.
		\end{equation}
		For all $x\in \L, y\in \ran P^\perp$ and $ n< 0,$  $\la U^nx,y\ra
		= \la  x, U^{-n}y\ra=0$  implies that		\begin{equation}\label{eq:n<0}
			\oplus_{n< 0}U^n(\L)\subseteq \ran P.
		\end{equation}
		It is clear from inclusions \eqref{eq:n>0} and \eqref{eq:n<0} that $P$ also reduces $\E_1.$
		Now $U|_{\E_1}=W_\L,$  the bilateral shift on $\E_1$ with the  wandering subspace $\L$ and  $P|_{\E_1}=P_{\E_1^-},$ the projection on $\E_1^-:=\oplus_{n< 0}U^n(\L)$ in $\E_1.$
		
		Let $\E=\ker V^*= \E_1\oplus \E_2.$ In this decomposition we have, \[U=\begin{pmatrix} W_\L &0 \\ 0 & U|_{\E_2}\end{pmatrix} \text{ and }  P=\begin{pmatrix}P_{\E_1^-}&0 \\ 0 & P|_{\E_2} \end{pmatrix}.\]
		Let $\Theta :\E_1\to l^2(\Z)\otimes \L$ be the unitary given by
		\begin{equation*}
			\Theta(U^nx)=e_n\otimes x \text{ for } x\in \L, n\in \Z.
		\end{equation*}
		Then, with this identification $(U,P)$ is jointly unitarily  equivalent to
		$$(\begin{pmatrix} \omega \otimes I_\L &0 \\ 0 & U|_{\E_2}\end{pmatrix} ,\begin{pmatrix} p_-\otimes I_\L&0 \\ 0 & P|_{\E_2} \end{pmatrix})$$
		where $\omega$ is the bilateral shift on $l^2(\Z)$ and $p_-$ is the projection in $l^2(\Z)$ onto $\spncl\{e_n:n<0\}.$
		Therefore,
		\begin{align}    \label{eq:negphi1z}
			\varphi_1(z)
			&=\begin{pmatrix}  \omega^*\otimes I_\L&0 \\ 0 & U^*|_{\E_2}\end{pmatrix} \left[ \begin{pmatrix}  p_-^\perp\otimes I_\L &0 \\ 0 & P^\perp|_{\E_2}\end{pmatrix}+z\begin{pmatrix} p_- \otimes I_\L&0 \\ 0 & P|_{\E_2}\end{pmatrix}\right]\nonumber\\
			&=\begin{pmatrix} \omega^*(p_-^\perp+zp_-)\otimes I_\L &0 \\ 0 & U^*(P^\perp+zP)|_{\E_2} \end{pmatrix}=\begin{pmatrix} \psi_1(z)\otimes I_\L&0 \\ 0 & \varphi_1(z)|_{\E_2} \end{pmatrix},
		\end{align}
		and
		\begin{align}   \label{eq:negphi2z}
			\varphi_2(z)
			&=\left[ \begin{pmatrix}  p_- \otimes I_\L&0 \\ 0 &P|_{\E_2}\end{pmatrix}+z\begin{pmatrix}  p_-^\perp\otimes I_\L &0 \\ 0 & P^\perp|_{\E_2} \end{pmatrix}\right]\begin{pmatrix} \omega \otimes I_\L &0 \\ 0 & U|_{\E_2} \end{pmatrix}\nonumber\\
			&=\begin{pmatrix} (p_-+zp_-^\perp)\omega \otimes I_\L&0 \\ 0 & (P+zP^\perp)U|_{\E_2} \end{pmatrix}=\begin{pmatrix} \psi_2(z)\otimes I_\L&0 \\ 0 & \varphi_2(z)|_{\E_2} \end{pmatrix},
		\end{align}
		where $p_-^\perp$ denotes the projection $I_{l^2(\Z)}-p_-$ in $l^2(\Z). $ This in particular says that $\E_1$ reduces $\varphi_i(z)$ for $z\in \D,i=1,2.$ Therefore $H^2_\D(\E_1)$ reduces $M_{\varphi_i},$ for $ i=1,2,$ and
		\begin{equation}
			M_{\varphi_i}=\begin{pmatrix}  M_{\psi_i\otimes I_\L }&0\\0& M_{{\varphi_i|_{\E_2}}} \end{pmatrix}
			=\begin{pmatrix} M_{\psi_i}\otimes I_\L &0\\0& M_{\varphi_i}|_{H^2_\D(\E_2)} \end{pmatrix}, \ i=1,2,
		\end{equation}
		where
		$$\psi_1(z)=\omega^*(p_-^\perp+zp_-),\  \psi_2(z)= (p_-+zp_-^\perp)\omega$$
		and $\varphi_i|_{\E_2}(z):=\varphi_i(z)|_{\E_2}\in\B(\E_2),$ for $z\in \D$ and $i=1,2.$
		
		
		It remains to prove that $C({M_{\varphi_1}}|_{H^2_{\D}(\E_2)},{M_{\varphi_2}}|_{H^2_{\D}(\E_2)})=0$. To that end, note that by Lemma \ref{lem:zeroequiv}, it is enough to show that $U(\ran P|_{\E_2})=\ran P|_{\E_2}.$ Notice from \eqref{eq:n>0} and \eqref{eq:n<0} that $\ran P|_{\E_2}=\ran P \ominus (\oplus_{n< 0}U^n(\L)).$ Hence
		\begin{align*}
			U(\ran P|_{\E_2})&=U(\ran P \ominus (\oplus_{n< 0}U^n(\L)))\\
			&=U(\ran P)\ominus (\oplus_{n< 0}U^{n+1}(\L))\\
			&=(U(\ran P^\perp))^\perp \ominus (\L\oplus_{n< 0}U^n(\L))\\
			&=(\ran P \oplus \L )\ominus (\L\oplus_{n< 0}U^n(\L))\quad (\because \ran P^\perp=\L\oplus U(\ran P^\perp))\\
			&=\ran P \ominus (\oplus_{n< 0}U^n(\L))=\ran P|_{\E_2}.
		\end{align*}
		It remains only to prove \eqref{eq:incpsi}, and it follows directly from \eqref{eq:mphimpsi} by Lemma \ref{ReducingSpectrum}.
	\end{proof}

	\begin{lemma}\label{lem:eq-eg}
		The pair $(M_{\psi_1},M_{\psi_2})$ is jointly unitarily equivalent to $(\tau_1,\tau_2).$ In particular, $$\sigma(M_{\psi_1},M_{\psi_2})=\overline{\D^2}$$ and for every point in $\D^2,$ the non-singularity breaks at stage 2.   
	\end{lemma}
	
	\begin{proof}


		Define the unitary $\Lambda:H^2_{\D^2}\to H^2_\D(l^2(\Z))$ by
		\begin{equation}
			\Lambda\left(\sum_{m,n=0}^{\infty}a_{m,n}z_1^mz_2^n\right)=\sum_{k=0}^\infty\left(\sum_{m= 0}^\infty a_{m+k,k}e_m+\sum_{m= 1}^\infty a_{k,m+k}e_{-m}\right)z^k.
		\end{equation}
		Then, $\Lambda \tau_i\Lambda^*=M_{\psi_i}$ for $i=1,2.$ That is, $(M_{\psi_1},M_{\psi_2})$ is jointly unitarily equivalent to $(\tau_1,\tau_2).$	In particular,
		\begin{equation}
			\sigma(M_{\psi_1},M_{\psi_2})=\sigma(\tau_1,\tau_2)=\overline{\D^2}
		\end{equation}
		and for every point in $\D^2,$ the non-singularity breaks at stage 2 by Proposition \ref{pro:eg}.
	\end{proof}
	We shall use the following lemma proved in \cite{GG} and \cite{Sandip}.
	\begin{lemma}\label{jointreducing}
		Let $(\E,P,U)$ be a BCL triple. Then the pair $(M_{\varphi_1},M_{\varphi_2})$ has a non-trivial joint reducing subspace if and only if $(P,U)$ has a non-trivial joint reducing subspace.
	\end{lemma}

	\begin{lemma} \label{Lem:Irreducible T1T2}
		The pair $(\tau_1,\tau_2)$ does not have any non-trivial joint reducing subspace.
	\end{lemma}
	
	\begin{proof}
		By Lemma \ref{jointreducing}, $(\tau_1,\tau_2)$ has a non-trivial joint reducing subspace if and only if $(p_-,\omega)$  has a non-trivial joint reducing subspace in $l^2(\Z).$ Let $\E_{0}\ne \{0\}$ be a joint reducing subspace for $(p_-,\omega).$ Let $f=\sum_{n\in \Z} a_ne_n\in \E_{0},$ $a_{n_0}\ne 0 $ for some $n_0\in \Z.$ Now
		\begin{equation*}
			{\omega^*}^{n_0}(f)=\sum_{n\in \Z} a_ne_{n-n_0},\ 	p_-{\omega^*}^{n_0}(f)=\sum_{n<n_0} a_ne_{n-n_0}\in \E_{0}.
		\end{equation*}
		Hence $\sum_{n\ge n_0} a_ne_{n-n_0}\in \E_{0}$ and
		\begin{equation*}
			p_-\omega^*(\sum_{n\ge n_0} a_ne_{n-n_0})=a_{n_0}e_{-1}\in \E_{0}.
		\end{equation*}
		Since $\E_{0}$ is reducing for $\omega$ and $e_{-1}\in \E_{0}$ we have $ \E_{0}=l^2(\Z).$
	\end{proof}
	
	The following theorem on the structure and the joint spectrum of commuting pair of isometries with negative defect follows directly from Theorem \ref{thm:BCL}, Theorem \ref{thm:NegStructure}, Lemma \ref{ReducingSpectrum} and Lemma \ref{lem:eq-eg} and shows that when the defect of $(V_1,V_2)$ is negative, then apart from a reduced part which has defect zero, $(V_1,V_2)$ is the fundamental isometric pair with negative defect, albeit with a higher multiplicity.
	
	\begin{thm}\label{StructureNegative}
		Let $(V_1,V_2)$ be a pair of commuting isometries such that $C(V_1,V_2)\le 0$ and $C(V_1,V_2)\ne 0.$  Then, there is a non-trivial subspace  $\L\subsetneq \ker V^*$ such that,  up to unitary equivalence, $\H=\H_0\oplus \H_0^\perp,$ where $\H_0= H^2_{\D^2}\otimes\L$ and in this decomposition 	
		\begin{equation*}\label{eq:phipsi}  V_i
			= \begin{pmatrix}  \tau_i\otimes I_\L &0\\0& V_i|_{\H_0^\perp} \end{pmatrix}, i=1,2 \text{ and } C(V_1|_{\H_0^\perp},V_2|_{\H_0^\perp})=0,
		\end{equation*}
		where the  dimension of $\L$ is same as the dimension of the range of $C(V_1,V_2)$ and  $(\tau_1, \tau_2)$ is the fundamental isometric pair with negative defect.  Moreover, $$\sigma(V_1,V_2)=\overline{\D^2}$$ and for every point in $\D^2,$ the non-singularity breaks at stage 2.
	\end{thm}

The pair $(\tau_1, \tau_2)$ serves as a model in many ways.

	\begin{thm}
	A pair of commuting isometries $(V_1, V_2)$ is a modified bi-shift if and only if $(V_1, V_2)$ is jointly unitarily equivalent to $(\tau_1 \otimes I_\L , \tau_2 \otimes I_\L)$ for some Hilbert space $\L$.
\end{thm}

	\begin{proof}
	Let $(V_1,V_2)=(\tau_1 \otimes I_\L , \tau_2 \otimes I_\L).$ Then
	$$\ker V_1^*=\spncl \{z_2^2,z_2^3,\dots\}\otimes \L,\  \ker V_2^*=\spncl\{1,z_1,z_1^2,\dots\}\otimes \L.$$
	Clearly $V_1^*|_{\ker V_2^*}, V_2^*|_{\ker V_1^*}$ and $V_1V_2$ are shifts. Therefore, by \cite[Prop. 3.8]{Popovici} $(V_1,V_2)$ is a modified bi-shift.
	
	Suppose $(V_1,V_2)$ is a modified bi-shift. If $(\E,P,U)$ is the BCL triple for $(V_1,V_2),$ then by \cite[Thm. 2.4]{GG}, $U(\ran P^\perp)\subsetneq \ran P^\perp.$ Hence by Lemma \ref{lem:UP}, $C(V_1,V_2)\le 0$ and $C(V_1,V_2)\ne 0.$ Therefore, by Theorem \ref{StructureNegative}, up to unitary equivalence, $V_1$ and $V_2$ are
	$$V_1=\begin{pmatrix}
		\tau_1\otimes I_\L &0& 0\\0&M_z\otimes I_{\E_1}& 0\\0&0& I_{H^2_\D}\otimes U_2
	\end{pmatrix}, V_2=\begin{pmatrix}
		\tau_2\otimes I_\L &0& 0\\0&I_{H^2_\D}\otimes U_1& 0\\0&0& M_z\otimes I_{\E_2}
	\end{pmatrix}.$$
	By \cite[Prop. 3.8]{Popovici}  we have $\E_1=\E_2=\{0\}.$ Hence $(V_1,V_2)$ is jointly unitarily equivalent to $(\tau_1\otimes I_\L,\tau_2\otimes I_\L).$
\end{proof}

\begin{remark} \label{rem:Neg}
	$(l^2(\Z)\otimes\L, p_-\otimes I_\L,\omega\otimes I_\L)$ is the BCL triple for a modified bi-shift $(V_1,V_2),$
	where $\dim \L=\dim (\ran C(V_1,V_2)). $ (See also \cite[Thm. 2.4]{GG}).
\end{remark}
 Some other descriptions of  modified bi-shift are given in  \cite{GG}, by Gaspar and Gaspar. From the results in \cite{GG}, it is clear that a dual doubly commuting pair is a direct sum of a zero defect part and a modified bi-shift. (caution: the $W$ in \cite{GG} is $U^*$ for us).

	In the rest of this section we shall see the relation between the joint spectrum of a pair of commuting isometries $(V_1,V_2)$ with  defect negative and non-zero,  and the joint spectrum of $(\varphi_1(z),\varphi_2(z)),z\in \D,$ where $\varphi_i$'s are the multipliers given in the BCL representation of $(V_1,V_2).$ Indeed, we prove that  \begin{equation}\label{eq:sp}
		\sigma(V_1,V_2)=\sigma(M_{\varphi_1},M_{\varphi_2})=\overline{\cup_{z\in \D}\sigma(\varphi_1(z),\varphi_2(z))}=\overline{\D^2}.
	\end{equation}
	
	\begin{lemma}\label{lem:negpsi}	Let the operator valued functions $\psi_1$ and $\psi_2$ be as in Theorem \ref{thm:NegStructure}. Then,
		For $z\in \D,$
		$ (\lambda_1,\lambda_2)\in\D^2$ is a joint eigenvalue for $(\psi_1(z),\psi_2(z))$ if  $\lambda_1\lambda_2=z.$  Also
		\begin{equation}
			\sigma(\psi_1(z),\psi_2(z))=\{(\lambda_1,\lambda_2)\in\overline{\D^2}:\lambda_1\lambda_2=z\},
		\end{equation}
		and
		\begin{equation}
			\label{spec clo}\overline{\cup_{z\in \D}\sigma(\psi_1(z),\psi_2(z))}
			=\overline{\D^2}=\sigma(M_{\psi_1},M_{\psi_2}).
		\end{equation}
	\end{lemma}
	
	\begin{proof}
		Let $z\in\D\setminus\{0\}$ and $\lambda_1,\lambda_2\in \D$ be such that $\lambda_1\lambda_2=z.$ Consider
		$$x_{\lambda_1,\lambda_2}=\sum_{n= 0}^\infty\lambda_1^ne_n+\frac{1}{\lambda_1}\sum_{n= 1}^\infty\lambda_2^{n-1}e_{-n}.$$
		Note that $$\psi_i(z)x_{\lambda_1,\lambda_2}=\lambda_ix_{\lambda_1,\lambda_2}, i=1,2.$$
		That is, $(\lambda_1,\lambda_2)$ is a joint eigenvalue for $(\psi_1(z), \psi_2(z))$ with eigenvector $x_{\lambda_1,\lambda_2},$ where $z\ne 0$ and $\lambda_1,\lambda_2\in \D$ are such that $\lambda_1\lambda_2=z.$ Therefore,
		\begin{equation*}
			\sigma(\psi_1(z),\psi_2(z))\supseteq\{(\lambda_1,\lambda_2)\in\overline{\D^2}:\lambda_1\lambda_2=z\} \text{ for } z\ne 0.
		\end{equation*}
		
		As $\psi_1(z)$ and $\psi_2(z)$ are commuting contractions, $\sigma(\psi_1(z),\psi_2(z))\subseteq \overline{\D^2}.$ Since  $\psi_1(z)\psi_2(z)=zI,$ by spectral mapping theorem,  we have the other inclusion. Hence
		\begin{equation}\label{eq:znot0}
			\sigma(\psi_1(z),\psi_2(z))=\{(\lambda_1,\lambda_2)\in\overline{\D^2}:\lambda_1\lambda_2=z\} \text{ for } z\ne 0.
		\end{equation}
		For the case $z=0,$ consider for $\lambda_1,\lambda_2\in\D\setminus\{0\},$
		$$x_{\lambda_1,0}=\sum_{n= 0}^\infty\lambda_1^ne_n+\frac{1}{\lambda_1}e_{-1},\ x_{0,\lambda_2}=\sum_{n= 1}^\infty\lambda_2^{n-1}e_{-n}\text{ and } x_{0,0}=e_{-1} .$$
		Then $(\lambda_1,0),(0,\lambda_2)$ and $(0,0)$ are joint eigenvalues for $(\psi_1(0),\psi_2(0))$ with the joint eigenvectors $x_{\lambda_1,0},x_{0,\lambda_2}$ and $x_{0,0}$ respectively.  Therefore, by the similar reasoning as in the case of $z\ne 0,$ we get
		\begin{equation}\label{eq:z=0}
			\sigma(\psi_1(0),\psi_2(0))=\Dc\times \{0\}\cup \{0\}\times \Dc=\{(\lambda_1,\lambda_2)\in\overline{\D^2}:\lambda_1\lambda_2=0\}.
		\end{equation}
		Hence, by \eqref{eq:znot0}, \eqref{eq:z=0} and Lemma \ref{lem:eq-eg}, we get
		\eqref{spec clo}.
		
		Indeed, one can show that the joint eigen spaces are one dimensional, for all the joint eigenvalues of $(\psi_1(z),\psi_2(z)), z\in \D.$
	\end{proof}
	
	\begin{thm}\label{thm:negspec}
		Let $(V_1,V_2)$ be a pair of commuting isometries such that $C(V_1,V_2)\le 0$ and $C(V_1,V_2)\ne 0.$ Let $(\E,P,U)$ be the BCL triple for $(V_1,V_2).$  Then for $z\in \D$, $ (\lambda_1,\lambda_2)\in\D^2$ is a joint eigenvalue for $(\varphi_1(z),\varphi_2(z))$ if  $\lambda_1\lambda_2=z,$
		\begin{equation*}
			\sigma(\varphi_1(z),\varphi_2(z))=\{(\lambda_1,\lambda_2)\in\overline{\D^2}:\lambda_1\lambda_2=z\},
		\end{equation*}
		and
		\begin{equation*}
			\sigma(V_1,V_2)=\sigma(M_{\varphi_1},M_{\varphi_2})=\overline{\cup_{z\in \D}\sigma(\varphi_1(z),\varphi_2(z))}=\overline{\D^2}.
		\end{equation*}
		Moreover, the non-singularity breaks at stage 2.
	\end{thm}
	
	\begin{proof}
		As $\varphi_1(z)$ and $\varphi_2(z)$ are commuting contractions, $\sigma(\varphi_1(z),\varphi_2(z))\subseteq\overline{\D^2}.$ Since  $\varphi_1(z)\varphi_2(z)=zI,$ by spectral mapping theorem, $$\sigma(\varphi_1(z),\varphi_2(z))\subseteq\{(\lambda_1,\lambda_2)\in\overline{\D^2}:\lambda_1\lambda_2=z\}.$$ From \eqref{eq:negphi1z} and \eqref{eq:negphi2z} note that:
		\begin{equation*}
			\varphi_i(z)=\begin{pNiceMatrix}[first-row,last-col,nullify-dots]
				l^2(\Z)\otimes\L & \E_2\\
				\psi_i(z)\otimes I_\L&0 & l^2(\Z)\otimes\L\\ 0 & \varphi_i(z)|_{\E_2}&\E_2	\end{pNiceMatrix}, \ i=1,2.
		\end{equation*}
		Now the proof follows from Lemma \ref{lem:negpsi}, Lemma \ref{ReducingSpectrum} and Lemma \ref{lem:eq-eg}.
	\end{proof}
	

	\section{The Positive Defect Case} \label{PositiveDefectSection}
   By Lemma \ref{lem:SLOBCL} (see also \cite{Yang,Sarkar}), $C(V_1,V_2)\ge 0$ if and only if $(V_1,V_2)$ is doubly commuting. The structure of doubly commuting pair of isometries are well understood in the literature; see \cite{Slocinski,GG,Popovici I}. We rephrase and shine some of the existing results using the defect operator. We also study the joint spectrum in detail. 	
	
	
	
	
	\subsection{The prototypical example}
	
	\begin{definition}
		The fundamental isometric pair of positive defect is the bi-shift $(M_{z_1}, M_{z_2}),$ the pair of multiplication by the coordinate functions on $H^2_{\D^2}$. \end{definition}
	It is folklore that the $M_{z_i}$ are isometries and the defect is positive because of the following lemma. It will be clear from Theorem \ref{StructurePositive} why we call it $fundamental$.
	\begin{lemma} \label{Defect_Hardy}
		The defect operator $C(M_{z_1}, M_{z_2})$ is the projection onto the one dimensional space of constant functions.
	\end{lemma}
	\begin{proof}
		The proof is a straightforward computation.
	\end{proof}
	
	If $(V_1,V_2)$ is a bi-shift on $\H$ with the wandering subspace $\mathcal R,$ then $\Lambda:\H\to H^2_{\D^2}\otimes \mathcal R$ given by $\Lambda (V_1^{n_1}V_2^{n_2}x)= z_1^{n_1}z_2^{n_2}\otimes x, x\in \mathcal{R}$ is a unitary and $\Lambda V_i \Lambda^*= M_{z_i}\otimes I_{\mathcal R}$ for $i=1,2.$ Also, by the above lemma $\dim (\ran C(V_1,V_2))=\dim \mathcal R.$ The following lemma is well known.
	
	\begin{lemma}\label{lem:concrete:eg}
		The joint spectrum $\sigma(M_{z_1}, M_{z_2})$ is the whole bidisc $\overline{\mathbb D^2}$. Indeed, every point in the open bidisc is a joint eigenvalue for $(M_{z_1}^*,M_{z_2}^*).$ In particular, for every $(w_1, w_2) \in \D^2$, the non-singularity of $K(M_{z_1}-w_1I,M_{z_2}-w_2I)$ is broken at the third stage. 
	\end{lemma}

	\begin{remark}\label{NrestrictedP}
		Recall the pair $(\tau_1, \tau_2)$ defined in Subsection \ref{subsec:neg:eq}.	One can observe that  $\{\tau_1^m\tau_2^n(1):m,n\ge0\}$ is an orthonormal subset of $H^2_{\D^2}.$ Let $\mathcal M=\spncl\{\tau_1^m\tau_2^n(1):m,n\ge 0\}\subset H^2_{\D^2}.$  Clearly it is a joint invariant subspace (but not reducing) for  $(\tau_1,\tau_2).$ Identify $\mathcal{M}$ and $H^2_{\D^2}$ via $\tau_1^m\tau_2^n(1)\mapsto z_1^mz_2^n.$ Clearly, the pair of isometries $(\tau_1|_\M,\tau_2|_\M)$ gets identified with $(M_{z_1},M_{z_2})$ on $H^2_{\D^2}$. Since $C(M_{z_1},M_{z_2})\ge 0$ and $C(M_{z_1},M_{z_2})\neq 0$, the same is true for the pair $(\tau_1|_\M,\tau_2|_\M)$. Thus, a pair of commuting isometries with negative defect, when restricted to a joint invariant subspace, can have positive defect.
	\end{remark}
	
	\subsection{The general case of the positive defect operator}
	We start with a characterization available in \cite{Yang} and \cite{Sarkar}. 
		
	\begin{lemma}\label{lem:PCI}
		Let $(V_1,V_2)$ be a pair of commuting isometries on a Hilbert space $\H$. Then the following are equivalent:
		\begin{itemize}
			\item[(a)] $C(V_1,V_2)\geq 0$ and $C(V_1,V_2)\neq 0.$
			\item[(b)] $V_2(\ker V_1^*)\subsetneq \ker V_1^*.$
			\item[(c)] $V_1(\ker V_2^*)\subsetneq \ker V_2^*.$
			\item[(d)] The fringe operators are isometries and not unitaries.
			\item[(e)] $(V_1,V_2)$ is doubly commuting and $C(V_1,V_2)\neq 0.$
			\item[(f)] $C(V_1,V_2)$ is a non-zero projection.
			\item[(g)] If $(\E,P,U)$ is the BCL triple for $(V_1,V_2),$ then $U(\ran P)\subsetneq \ran P.$
		\end{itemize}
	\end{lemma}

 Let $\eta_1,\eta_2:\D\to \B(l^2(\Z))$ be the multipliers associated with the BCL triple $(l^2(\Z),p_{0+},\omega)$ where $p_{0+}$ is  the projection onto $\spncl\{e_n:n\ge0\}$ and $\omega$ is the bilateral shift on $l^2(\Z)$. 	Now, we obtain a structure theorem which has its own independent interest and is applied later to obtain Theorem \ref{thm:pos:sp:mult}.
	
	\begin{thm}\label{thm:PosStructure}
		Let $(V_1,V_2)$ be a pair of commuting isometries such that $C(V_1,V_2)$ is a non-zero positive operator. Let  $(\E,P,U)$ be the  BCL triple for $(V_1,V_2).$ Then, up to  unitary equivalence $$\E = (l^2(\Z)\otimes \L) \oplus \E_2$$
		for some non-trivial closed subspace $\L$ and a closed subspace $\E_2$ of $\E.$  Moreover
		\begin{equation}
			M_{\varphi_i}=\begin{pNiceMatrix}[first-row,last-col,nullify-dots]
				H^2_\D(l^2(\Z))\otimes \L& H^2_\D(\E_2)&\\
				M_{\eta_i}\otimes I_\L & 0&H^2_\D(l^2(\Z))\otimes \L \\ 0& {M_{\varphi_i}}|_{H^2_{\D}(\E_2)}& H^2_\D(\E_2)
			\end{pNiceMatrix}
		\end{equation}
		with $C({M_{\varphi_1}}|_{H^2_{\D}(\E_2)},{M_{\varphi_2}}|_{H^2_{\D}(\E_2)})=0.$ In particular,
		\begin{equation}
			\sigma(M_{\eta_1},M_{\eta_2})= \sigma(M_{\eta_1}\otimes I_\L, M_{\eta_2}\otimes I_\L)\subseteq\sigma(M_{\varphi_1},M_{\varphi_2})\subseteq \sigma(V_1,V_2).
		\end{equation}
	\end{thm}

	\begin{proof}
		We shall not give details of this proof because it follows the same line as the proof of Theorem \ref{thm:NegStructure}.
	\end{proof}
	
	\begin{lemma}\label{lem:pos:eta:Mzi}
		The pair $(M_{\eta_1},M_{\eta_2})$ is jointly unitarily equivalent to $(M_{z_1},M_{z_2}).$ In particular, $\sigma(M_{\eta_1},M_{\eta_2})=\overline{\D^2}.$
	\end{lemma}
	
	\begin{proof}
		
		
		Define the unitary $\Lambda:H^2_{\D^2}\to H^2_\D(l^2(\Z))$ by
		\begin{equation}
			\Lambda\left(\sum_{m,n=0}^{\infty}a_{m,n}z_1^mz_2^n\right)=\sum_{k=0}^\infty\left(\sum_{m= 0}^\infty a_{m+k,k}e_{-(m+1)}+\sum_{m= 1}^\infty a_{k,m+k}e_{m-1}\right)z^k.
		\end{equation}
		Then, $\Lambda M_{z_i}\Lambda^*=M_{\eta_i}$ for $i=1,2.$ That is, $(M_{\eta_1},M_{\eta_2})$ is jointly unitarily equivalent to $(M_{z_1},M_{z_2}).$
		
		In particular, $\sigma(M_{\eta_1},M_{\eta_2})=\sigma(M_{z_1},M_{z_2})=\overline{\D^2},$ by Lemma \ref{lem:concrete:eg}.
	\end{proof}

\begin{remark} \label{rem:Pos}
	$(l^2(\Z)\otimes \L, p_{0+}\otimes I_\L,\omega\otimes I_\L)$ is the BCL triple for the bi-shift with wandering subspace $\L.$
	(Compare this with the Remark \ref{rem:Neg}).
\end{remark}

	\begin{lemma}
		The pair $(M_{z_1},M_{z_2})$ does not have any non-trivial joint reducing subspace.
	\end{lemma}
The proof of the above lemma is similar to the proof of Lemma \ref{Lem:Irreducible T1T2}. Now the following structure theorem for $C(V_1,V_2)\ge 0$ follows from Theorem \ref{thm:BCL}, Theorem \ref{thm:PosStructure}, Lemma \ref{lem:pos:eta:Mzi}  and Lemma \ref{Defect_Hardy}. It shows the role played by the bi-shift of a possibly higher multiplicity and is a rephrasing of some results in \cite{Slocinski,GG} using the defect operator.
	\begin{thm} \label{StructurePositive}
		Let $(V_1, V_2)$ be a pair of commuting isometries on $\H$ with $C(V_1, V_2)\ge 0$ and $C(V_1, V_2)\ne 0$. Then, there is a non-trivial Hilbert space $\L \subsetneq \ker V^*$ such that, up to unitary equivalence,
		$\H=\H_0\oplus \H_0^\perp,$ where $\H_0= H^2_{\D^2}\otimes\L.$ In this decomposition
		$$ V_i = \left( \begin{array}{cc} M_{z_i}\otimes  I_\L& 0 \\
			0 & V_{i0}\end{array} \right)$$
		and the defect operator $C(V_{10}, V_{20})$ is zero.  Moreover, the dimension of $\L
		$ is the same as the dimension of the range of $C(V_1, V_2)$.
	\end{thm}
	

	%
	%
	
	
	\begin{thm} \label{JSPositive}
		Let $(V_1,V_2)$ be a pair of commuting isometries on a Hilbert space $\H$ with $C(V_1,V_2)\geq 0$ and $C(V_1,V_2)\neq 0.$ Then, $\sigma(V_1,V_2)=\overline{\D^2}$ and every point $(w_1,w_2)\in \D^2$ is a joint eigenvalue of $(V_1^*,V_2^*)$.   In particular, for every $(w_1, w_2) \in \D^2$, the non-singularity of $K(V_1-w_1I,V_2-w_2I)$ is broken at the third stage.
	\end{thm}
	
	\begin{proof} The proof follows from Theorem \ref{JSPositive} and Lemma \ref{lem:concrete:eg}.
	\end{proof}
	
	The following lemma is useful to see the relation between the joint spectrum of a pair of  commuting isometries $(V_1,V_2)$ with  defect positive and non-zero,  and the joint spectra $\sigma(\varphi_1(z),\varphi_2(z)),$ $ z\in \D,$ where $\varphi_i$'s are the multipliers given in the BCL theorem. 
	
	\begin{lemma}\label{lem:poseta}
		Let the operator valued functions $\eta_1$ and $\eta_2$ be as in Theorem \ref{thm:PosStructure}. Then, for $z\in \D,$	$ (\overline{\lambda_1},\overline{\lambda_2})\in\D^2$ is a joint eigenvalue for $(\eta_1(z)^*,\eta_2(z)^*)$ if  $\lambda_1\lambda_2=z.$ Also,
		\begin{equation}
			\sigma(\eta_1(z),\eta_2(z))=\{(\lambda_1,\lambda_2)\in\overline{\D^2}:\lambda_1\lambda_2=z\},
		\end{equation}
		and
		\begin{equation}
			\label{spec clos 2}\overline{\cup_{z\in \D}\sigma(\eta_1(z),\eta_2(z))}
			=\overline{\D^2}=\sigma(M_{\eta_1},M_{\eta_2}).
		\end{equation}
	\end{lemma}
	
	\begin{proof}
		Let $z\in\D\setminus\{0\}$ and $\lambda_1,\lambda_2\in \D$ be such that $\lambda_1\lambda_2=z.$ Consider
		$$x_{\lambda_1,\lambda_2}=\sum_{n= 0}^\infty\overline{\lambda_2}^{n}e_n+\frac{1}{\overline{\lambda_2}}\sum_{n= 1}^\infty\overline{\lambda_1}^{n-1}e_{-n}.$$
		Note that  $$\eta_i(z)^*x_{\lambda_1,\lambda_2}=\overline{\lambda}_ix_{\lambda_1,\lambda_2}, i=1,2.$$
		That is, $(\overline{\lambda_1},\overline{\lambda_2})$ is a joint eigenvalue for $(\eta_1(z)^*, \eta_2(z)^*)$ with eigenvector $x_{\lambda_1,\lambda_2},$ where $z\ne 0$ and $\lambda_1,\lambda_2\in \D$ are such that $\lambda_1\lambda_2=z.$ Therefore,
		\begin{equation*}
			\sigma(\eta_1(z),\eta_2(z))\supseteq\{(\lambda_1,\lambda_2)\in\overline{\D^2}:\lambda_1\lambda_2=z\} \text{ for } z\ne 0.
		\end{equation*}
		As $\eta_1(z)$ and $\eta_2(z)$ are commuting contractions, we have $\sigma(\eta_1(z),\eta_2(z))\subseteq \overline{\D^2}.$
		Since $\eta_1(z)\eta_2(z)=zI,$ by spectral mapping theorem, we have
		the other inclusion. Hence
		\begin{equation}\label{eq:znot0case}
			\sigma(\eta_1(z),\eta_2(z))=\{(\lambda_1,\lambda_2)\in\overline{\D^2}:\lambda_1\lambda_2=z\} \text{ for } z\ne 0.
		\end{equation}
		For the case $z=0,$ consider for $\lambda_1,\lambda_2\in\D\setminus\{0\},$
		$$x_{\lambda_1,0}=\sum_{n=1}^\infty{\overline{\lambda_1}}^{n-1}e_{-n},\ x_{0,\lambda_2}=\sum_{n=0}^\infty\overline{\lambda_2}^ne_{n}+\frac{1}{\overline{\lambda_2}}e_{-1}\text{ and } x_{0,0}=e_{-1} .$$
		Then $(\overline{\lambda_1},0),(0,\overline{\lambda_2})$ and $(0,0)$ are joint eigenvalues for $(\eta_1(0)^*,\eta_2(0)^*)$ with the joint eigenvectors $x_{\lambda_1,0},x_{0,\lambda_2}$ and $x_{0,0}$ respectively.
		Therefore, by the similar reasoning as in the case of $z\ne 0,$ we get
		\begin{equation}\label{eq:z=0case}
			\sigma(\eta_1(0),\eta_2(0))=\Dc\times \{0\}\cup \{0\}\times \Dc=\{(\lambda_1,\lambda_2)\in\overline{\D^2}:\lambda_1\lambda_2=0\}.
		\end{equation}
		Hence, by \eqref{eq:znot0case}, \eqref{eq:z=0case} and Lemma \ref{lem:pos:eta:Mzi}, we get \eqref{spec clos 2}.

		Indeed, one can show that the joint eigenspaces are one dimensional, for all the joint eigenvalues of $(\eta_1(z)^*,\eta_2(z)^*), z\in \D.$	
	\end{proof}
	
	\begin{thm}\label{thm:pos:sp:mult}
		Let $(V_1,V_2)$ be a pair of commuting isometries such that $C(V_1,V_2)\ge 0$ and $C(V_1,V_2)\ne 0.$ Let $(\E,P,U)$ be the BCL triple for $(V_1,V_2).$  Then, for $z\in \D,$ $ (\overline{\lambda_1},\overline{\lambda_2})\in\D^2$ is a joint eigenvalue for $(\varphi_1(z)^*,\varphi_2(z)^*)$ if  $\lambda_1\lambda_2=z,$   also
		\begin{equation}
			\sigma(\varphi_1(z),\varphi_2(z))=\{(\lambda_1,\lambda_2)\in\overline{\D^2}:\lambda_1\lambda_2=z\}.
		\end{equation}
		Moreover, every point in the open bidisc is a joint eigenvalue for $(M_{\varphi_1}^*,M_{\varphi_2}^*),$ and
		\begin{equation}
			\sigma(V_1,V_2)=\sigma(M_{\varphi_1},M_{\varphi_2})=\overline{\cup_{z\in \D}\sigma(\varphi_1(z),\varphi_2(z))}=\overline{\D^2}.
		\end{equation}
	\end{thm}
	
	\begin{proof}
		
		Since $\varphi_1(z)$ and $\varphi_2(z)$ are commuting contractions, $\sigma(\varphi_1(z),\varphi_2(z))\subseteq\overline{\D^2}.$ As $\varphi_1(z)\varphi_2(z)=zI,$  by spectral mapping theorem,  $$\sigma(\varphi_1(z),\varphi_2(z))\subseteq\{(\lambda_1,\lambda_2)\in\overline{\D^2}:\lambda_1\lambda_2=z\}.$$ In the same manner as in the proof of Theorem \ref{thm:NegStructure},	we get
		\begin{equation*}
			\varphi_i(z)=\begin{pNiceMatrix}[first-row,last-col,nullify-dots]
				l^2(\Z)\otimes \L & \E_2\\
				\eta_i(z)\otimes I_\L&0 & l^2(\Z)\otimes \L\\ 0 & \varphi_i(z)|_{\E_2}&\E_2	\end{pNiceMatrix}
			,\ i=1,2,
		\end{equation*}
		where $\L=\ran P\ominus U(\ran P)$ and $\E_2=\ker V^*\ominus ( l^2(\Z)\otimes\L).$ Now the proof follows from Lemma \ref{lem:poseta}, Lemma \ref{lem:pos:eta:Mzi} and Lemma \ref{lem:concrete:eg}.
	\end{proof}

	\section{Towards the general defect operator}
	This section deals with the cases of  $\ran V_1=\ran V_2$ and $\ran V_2\subsetneq \ran V_1.$ We provide  the characterization and study the joint spectrum for both the cases. At the end of this section we give an example for the unknown case of $\H_i\ne 0$ for all $i=1,2,3,4.$

	\subsection{Range  of $V_1$ equal  to the Range of $V_2$}
	

In this case the defect operator is the difference of two mutually orthogonal projections whose ranges together span the kernel of $V^*$. The structure in this case, known from \cite{Burdak}, is briefly recalled below because it is needed for deciphering the joint spectrum.
	

	\subsubsection{The prototypical family of examples}\label{ProtypicalOffDiag}

	Let $\L$ be a Hilbert space and $W$ be a unitary on $\L$. Consider the pair of commuting isometries $(M_z\otimes I, M_z\otimes W)$ on $H^2_{\D}\otimes \L$. Clearly $\ran (M_z\otimes I)=\ran (M_z\otimes W).$
	
	\begin{lemma}
		Let $\L$ and $W$ be as above. If $(V_1, V_2) = (M_z\otimes I, M_z\otimes W)$, then there are two projections $P$ and $Q$ such that
		\begin{enumerate}
			\item they are mutually orthogonal to each other,
			\item the dimensions of ranges of $P$ and $Q$ are same,
			\item the span of the ranges of $P$ and $Q$ is the kernel of $V^*$ and
			\item the defect operator $C(V_1, V_2)$ on $H^2_{\D}\otimes \L$ is $P - Q$.
		\end{enumerate}
	\end{lemma}
	
	\begin{proof}
		In keeping with the notation $E_0$ for the projection onto the one dimensional subspace of constants in $H^2_{\D}$, let us  denote by $E_1$ the projection onto the one dimensional space spanned by $z$. The kernel of $V^*$ is
		$$\ran(I \otimes I_\L - M_z^2(M_z^2)^* \otimes I_\L) = \ran((E_0 + E_1) \otimes I_\L)=\spncl\{1,z\}\otimes \L$$
		and the defect operator is
		\begin{align*}
			C(V_1, V_2) & = I \otimes I_\L - M_zM_z^* \otimes I_\L - M_zM_z^* \otimes I_\L + M_z^2(M_z^2)^* \otimes I_\L \\
			& = E_0 \otimes I_\L - M_zE_0M_z^* \otimes I_\L \\
			& = E_0 \otimes I_\L - E_1 \otimes I_\L. \end{align*}
		Thus, setting $P = E_0 \otimes I_\L $ and $Q = E_1 \otimes I_\L $, we are done.
	\end{proof}
	
	\begin{lemma}If $\L\ne \{0\},$ then the joint spectrum of $(M_{z}\otimes I_{\L}, M_{z}\otimes W)$ is $$\sigma(M_{z}\otimes I_{\L}, M_{z}\otimes W)=\{z(1,\alpha): z\in \Dc, \alpha \in \sigma (W)\}.$$
	\end{lemma}
	
	\begin{proof}
		This is a straightforward application of the polynomial spectral mapping theorem. Consider the polynomial  $f(x_1,x_2)=(x_1, x_1x_2).$ Then,
		\begin{align}\label{proto off}\sigma(M_{z}\otimes I_{\L}, M_{z}\otimes W)&=f(\sigma(M_{z}\otimes I_{\L}, I_{H^2_\D}\otimes W))\nonumber\\
			&=f(\sigma(M_{z})\times \sigma(W)) \nonumber\\
			&=\{z(1,\alpha): z\in \Dc, \alpha \in \sigma (W)\}.
		\end{align}
	\end{proof}
	
	\subsubsection{General Theory}
	
	\begin{thm}\label{char off}
		Let $(V_1,V_2)$ be a pair of commuting isometries on a Hilbert space $\H$. Then the following are equivalent:
		\begin{itemize}
			\item[(a)] $\ran V_1=\ran V_2.$
			\item[(b)] $V_1(\ker V_2^*)=V_2(\ker V_1^*).$
			\item[(c)]The defect operator $C(V_1,V_2)$ is a difference of two mutually orthogonal projections $Q_1, Q_2$ with $\ran Q_1\oplus \ran Q_2=\ker V^*.$
			\item[(d)] The fringe operators $F_1$ and $F_2$ are zero.
			\item[(e)] If $(\E,P,U)$ is the BCL triple for $(V_1,V_2),$ then $U(\ran P)=\ran P^\perp$ (or equivalently $U(\ran P^\perp)=\ran P$).
		\end{itemize}
	\end{thm}
	
	\begin{proof}
		
		\leavevmode

		$(a)\Rightarrow (b)$ This is immediate from the equality
		$$\ker V_1^*\oplus V_1(\ker V_2^*)= \ker V^*=\ker V_2^*\oplus V_2(\ker V_1^*)$$
		after we note that  $\ran V_1=\ran V_2$ implies that $\ker V_1^* = \ker V_2^*.$
		
		$(b)\Rightarrow(c)$ Set $P_1 = P_{\ker V_1^*}$ and $P_2 = P_{V_2(\ker V_1^*)}.$ By \eqref{eq:def:diffofproj}, $C(V_1,V_2)=P_1- P_2$. Moreover,
		\begin{align*}
			\ran P_1\oplus\ran P_2&=\ker V_1^*\oplus V_2(\ker V_1^*)\\
			&=\ker V_1^*\oplus V_1(\ker V_2^*) \quad \text{(by (b))}\\
			&=\ker V^*.
		\end{align*}


		$(c)\Rightarrow(d)$ It is immediate by noticing that
		$$Q_1=P_{\ker V_1^*}=P_{\ker V_2^*}\text{ and }Q_2=P_{V_2(\ker V_1^*)}=P_{V_1(\ker V_2^*)}$$ from the discussions in  Section \ref{sec:defopr}.
		
		$(c) \Rightarrow (a)$ Immediate from the last line above.
		
		$(d)\Rightarrow(c)$ $F_1=0$ and $F_2=0$ implies that $\ker V_1^*\perp V_2(\ker V_1^*)$  and  $\ker V_2^*\perp V_1(\ker V_2^*).$ Now from \eqref{eq:directsum}, we see that $V_1(\ker V_2^*)=V_2(\ker V_1^*).$ Set $P_1=P_{\ker V_1^*}$ and $P_2=P_{V_2(\ker V_1^*)},$ to obtain $(c).$
		
		Thus, we have shown the equivalences of $(a),\ (b),\ (c)$ and $(d).$
		
		
		$(a)\iff (e)$ If $(\E,P,U)$ is the BCL triple for $(V_1,V_2),$ then
		$$\ran V_1  =  \ran V_2$$
		if and only if
		$$\ker V_1^*  =  \ker V_2^*$$
		if and only if
		$$\ker M_{\varphi_1}^*  =  \ker M_{\varphi_2}^*$$
		if and only if
		$$\ran E_0\otimes U^*PU=\ran E_0\otimes P^\perp$$
		if and only if
		$$ U^*PU  =  P^\perp.$$
	\end{proof}
	
	The following result is a special case of Remark 4.2 in \cite{Burdak}, which shows that the pair $(M_z\otimes I,M_z\otimes W),$ studied in Subsection \ref{ProtypicalOffDiag}, is the prototypical pair in this case. 
	We give a proof which is different from the one in \cite{Burdak}.

	\begin{thm}\label{thm:OffDiagStruc}
		Let $(V_1,V_2)$ be a pair of commuting isometries on $\H$ satisfying any of the equivalent conditions in Theorem \ref{char off}. Then, there exist Hilbert spaces $\L$ and $\K$ such that up to  unitarily equivalence $\H=(H^2_{\D}\otimes \L) \oplus \K$ and in this decomposition,
		$$V_1=\begin{pmatrix}
			M_z\otimes I_{\L}& 0\\
			0&W_1
		\end{pmatrix}, \quad V_2= \begin{pmatrix}
			M_z\otimes W& 0\\
			0&W_2
		\end{pmatrix},$$
		for some unitary $W$ on $\L$ and commuting unitaries $W_1, W_2$ on $\K.$
	\end{thm}
	
	\begin{proof}
		Let $(\E,P,U)$ be the BCL triple for $(V_1,V_2).$ Let $\L=\ran P^\perp.$ By Theorem \ref{char off} $(e)$, we have $U(\L)={\L}^{\perp}$ and $U(\L^{\perp})={\L}.$ Let $U_1=U|_\L:\L \to \L^\perp$ and $U_2=U|_{\L^\perp}:\L^\perp \to \L.$  In the decomposition $\E={\L}^{\perp}\oplus \L,$ we have $$U=\begin{pmatrix}
			0 & U_1 \\
			U_2 & 0
		\end{pmatrix} \text{ and } P=\begin{pmatrix}
			I_{\L^{\perp}} & 0 \\
			
			0 & 0
		\end{pmatrix}.$$ Hence
		$\varphi_1(z)=\begin{pmatrix}
			0 & U_2^{*} \\
			z U_1^{*} & 0
		\end{pmatrix}$ and $\varphi_2(z)=\begin{pmatrix}
			0 &  U_1  \\
			z  U_2 & 0
		\end{pmatrix}$ for $z\in \D,$ and
		$$M_{\varphi_1}=\begin{pmatrix}
			0 & I\otimes U_2^* \\
			M_{z}\otimes U_1^* & 0
		\end{pmatrix} \text{ and }  M_{\varphi_2}=\begin{pmatrix}
			0 & I\otimes U_1 \\
			M_z\otimes U_2 & 0
		\end{pmatrix}.$$
		Let us define the unitary $$\Lambda: (H^2_\D\otimes {\L^\perp}) \oplus (H^2_\D\otimes \L)\rightarrow H^2_\D\otimes \L$$
		given by,
		$$\Lambda(\begin{pmatrix}
			\sum_{n=0}^{\infty} a_{n}z^{n} \\
			\sum_{n=0}^{\infty} b_{n}z^{n}
		\end{pmatrix}):=\sum_{n \text{ is even }} {(U_2U_1)}^{\frac{n}{2}}(b_{\frac{n}{2}})z^{n}+\sum_{n \text{ is odd }} {(U_2U_1)}^{\frac{n-1}{2}}U_2(a_{\frac{n-1}{2}})z^{n},$$ $a_n\in\L^\perp,b_n\in\L.$ One can see that $\Lambda M_{\varphi_1}{\Lambda}^{*}= M_{z}\otimes I_{\L}$ and $\Lambda M_{\varphi_2}{\Lambda}^{*}= M_{z}\otimes W,$ where $W=U_2U_1.$ This completes the proof by Theorem \ref{thm:BCL}.
	\end{proof}
	\begin{corollary}
		Let $(V_1,V_2)$ be a pair of commuting isometries satisfying any of the equivalent conditions in Theorem \ref{char off}. If $V_1V_2$ is pure, then both $V_1$ and $V_2$ are pure.
	\end{corollary}
	
	Since the BCL triple  is unique (up to unitary equivalence), we shall use the following convenient choice of BCL triple due to  A. Maji et al in \cite{Sarkar} for a pair of commuting isometries.
	
	\begin{thm}\label{thm:Sarkar}
		Let $(V_1,V_2)$ be a pair of commuting isometries on a Hilbert space $\H$. Then, $(\ker V^*,P,U_0)$ is the BCL triple for $(V_1,V_2),$ where $P\in \B(\ker V^{*})$ is the orthogonal projection onto $V_2(\ker V_1^*)$ and
		$$U_0=\begin{pmatrix}
			V_2|_{\ker V_1^*} & 0 \\
			0 & V_1^*|_{V_1(\ker V_2^*)}
		\end{pmatrix}: \begin{array}{ccc}
			\ker V_1^*&     &V_2(\ker V_1^*)  \\
			\oplus &\to   &\oplus\\
			V_1(\ker V_2^*)& &\ker V_2^*
		\end{array}$$ is a unitary operator on $\ker V^*.$
	\end{thm}
	
	\begin{remark}\label{Rem:off}
		In Theorem \ref{thm:OffDiagStruc}, we can write the $W$ and $\L$ explicitly in terms of $V_1$ and $V_2$ as follows:
		
		By Theorem \ref{thm:Sarkar} and Theorem \ref{char off},  $\L=\ran P^\perp=\ker V_1^*,$	$U_1:\ker V_1^*\to V_1(\ker V_1^*)$ given by $U_1=V_2|_{\ker V_1^*}$ and $U_2:V_1(\ker V_1^*)\to \ker V_1^*$ given by $U_2=V_1^*|_{V_1(\ker V_1^*)}.$ Therefore $W=V_1^*|_{V_1(\ker V_1^*)}V_2|_{\ker V_1^*}.$
	\end{remark}
	\subsubsection{Joint spectrum}
	If $(V_1,V_2)$ is pure and satisfying any of the equivalent conditions in Theorem \ref{char off}, then by Theorem \ref{thm:OffDiagStruc},  Remark \ref{Rem:off} and by subsection \ref{proto off}, we have
	\begin{equation}
		\sigma(V_1,V_2)
		=\{z(1,e^{i\theta}): z\in \overline{ \D}, \ e^{i\theta}\in \sigma(U_2U_1)\},
	\end{equation}
	where $U_1=V_2|_{\ker V_1^*},U_2=V_1^*|_{V_1(\ker V_1^*)}.$
	
	The final theorem of this section tells us the nature of elements in the joint spectrum and the relation between the joint spectrum of the commuting isometries satisfying any of the equivalent conditions in Theorem \ref{char off} and the joint spectra of the associated multipliers at every point of $\D.$

	\begin{thm}\label{offspecmulti}	Let $(V_1,V_2)$ be a pair of commuting isometries on $\H$ satisfying any of the equivalent conditions in Theorem \ref{char off} and $\ker V^*\ne \{0\}.$ Let $(\E,P,U)$ be the BCL triple for $(V_1,V_2).$    Let
		$U_1:\ran P^\perp \to \ran P \text{ be the unitary given by } U_1(x)=U(x), x\in \ran P^\perp$ and $U_2:\ran P \to \ran P^\perp \text{ be the unitary given by } U_2(y)=U(y), y\in \ran P.$
		Then
		\begin{itemize}
			\item $\sigma(\varphi_1(z),\varphi_2(z))=\{\pm \sqrt{z} (e^{-i \frac{\theta}{2}} , e^{i \frac{\theta}{2}}) : e^{i \theta} \in \sigma(U_1U_2) \}.$
		Here for every point in the joint spectrum, the non-singularity breaks at stage 3.
		\item		$\sigma(M_{\varphi_1},M_{\varphi_2})=\overline{\cup_{z\in \D}\sigma(\varphi_1(z),\varphi_2(z))}= \{z(1,e^{i\theta}): z\in \overline{ \D}, \ e^{i\theta}\in \sigma(U_1U_2)\}.$		
		Here, for every point $(z_1,z_2)$ in the set  $\{z(1,e^{i\theta}): z\in \D, \ e^{i\theta}\in \sigma(U_1U_2)\}$,
		the non-singularity in the Koszul complex $K(M_{\varphi_1}-z_1I,M_{\varphi_2}-z_2I)$   breaks at stage 3.		
	\end{itemize}	
	\end{thm}

	
	\begin{proof}
		Let $\E_1=\ran P$ and $\E_2=\ran P^\perp.$
		We have by Theorem \ref{char off} that $U(\E_1)=\E_2$ and  $U(\E_2)=\E_1.$ Hence in the decomposition $\E=\E_1\oplus \E_2,$ we have
		$\varphi_1(z)=\begin{pmatrix}
			0 & {U_2^{*}} \\
			z {U_1^{*}} & 0
		\end{pmatrix}$ and $\varphi_2(z)=\begin{pmatrix}
			0 & U_1 \\
			z U_2 & 0
		\end{pmatrix},$ for $z\in \D,$ where $U_1$ and $U_2$ are as in the statement.
		
		Now, we shall show that: for $z\in \D,$
		\begin{align}
			\sigma(\varphi_1(z))&=\{w\in \D: w^{2}=\overline{\alpha}z\text{ for some }  \alpha \in \sigma(U_1U_2)\},\label{eq:sp:phi1}\\ \sigma(\varphi_2(z))&=\{w\in \D: w^{2}={\alpha}z\text{ for  some } \alpha \in \sigma(U_1U_2)\}.\label{eq:sp:phi2}
		\end{align}
		Clearly  \eqref{eq:sp:phi1} and \eqref{eq:sp:phi2} hold for $z=0$. So assume that $z\neq 0$.
		
		\vspace*{3mm}
		
		\noindent {\bf Assertion 1.} If $\lambda$ is an eigenvalue of $\varphi_1(z)$ then $\frac{\overline{\lambda^{2}}}{\overline{z}}$ is an eigenvalue of $U_1U_2$.  In particular, $\lambda \in \D$.
		
		\vspace*{3mm}
		
		\noindent{\bf Proof of Assertion 1.} If $\lambda$ is an eigenvalue of $\varphi_1(z),$  there exists a non-zero vector $h_1\oplus h_2
		\in \E_1\oplus \E_2$ such that $\varphi_1(z)\begin{pmatrix}
			h_1 \\
			h_2
		\end{pmatrix}=\lambda \begin{pmatrix}
			h_1 \\
			h_2
		\end{pmatrix}$ which implies ${U_2}^{*}(h_2)=\lambda h_1$ and $\frac{z}{\lambda} {U_1^{*}}(h_1)= h_2.$ Hence ${(U_1U_2)^{*}}h_1=\frac{{\lambda^{2}}}{z}h_1$. Notice that $h_1\neq 0$, otherwise ${U_2}^{*}(h_2)=\lambda h_1$ implies $h_2=0$.
		
		\vspace*{3mm}
		
		\noindent {\bf Assertion 2.} For $\lambda \in \mathbb{C}$, $(\varphi_1(z)-\lambda I_{\E})$ is not onto if and only if $({(U_1U_2)^{*}}-\frac{{\lambda^{2}}}{z} I_{\E_1})$ is not onto. In particular, if $(\varphi_1(z)-\lambda I_{\E})$ is not onto then  $\lambda \in \D$.
		Also
		\begin{equation}\label{eq:ran:phi1}
			\ran(\varphi_1(z)-\lambda I)=\{\begin{pmatrix}
				h_1 \\
				h_2
			\end{pmatrix}\in \E_1\oplus \E_2: \lambda h_1+ U_2^* h_2\in  \ran ((U_1U_2)^*-\frac{\lambda^2}{z} I) \}.
		\end{equation}
		
		\vspace*{3mm}
		
		\noindent {\bf Proof of Assertion 2.}	Let $\begin{pmatrix}
			h_1 \\
			h_2
		\end{pmatrix}\in \E_1\oplus \E_2.$ We have $(\varphi_1(z)-\lambda I)\begin{pmatrix}
			h_1 \\
			h_2
		\end{pmatrix}=\begin{pmatrix}
			-\lambda h_1+{U_2}^{*} h_2 \\
			z {U_1}^{*} h_1- \lambda h_2
		\end{pmatrix}.$
		Let $k=\begin{pmatrix}
			k_1 \\
			k_2
		\end{pmatrix}\in \E_1\oplus \E_2.$  Now $(\varphi_1(z)-\lambda I)\begin{pmatrix}
			h_1 \\
			h_2
		\end{pmatrix}=\begin{pmatrix}
			k_1 \\
			k_2
		\end{pmatrix}$ if and only if $-\lambda h_1 +{U_2}^{*} (h_2)=k_1$ and
		$$({(U_1U_2)^{*}}-\frac{{\lambda^{2}}}{z} I_{\E_1})(h_1)=\frac{\lambda k_1+{U_2}^{*}(k_2)}{z}.$$
		This proves Assertion 2.

		Using Assertions 1 and 2 and the fact that for any bounded normal operator $T$ if $\alpha \in \sigma(T)$ then $T- \alpha I$ is not onto, we get \eqref{eq:sp:phi1} for $z\ne 0$. In a similar way, one can show \eqref{eq:sp:phi2} and for $z\neq 0$
		\begin{equation}\label{eq:ran:phi2}
			\ran(\varphi_2(z)-\lambda I)=\{\begin{pmatrix}
				k_1 \\
				k_2
			\end{pmatrix}\in \E_1\oplus \E_2: \lambda k_1+ U_1 k_2\in \ \ran ((U_1U_2)-\frac{{\lambda^{2}}}{z} I) \}.	
		\end{equation}
		For $z\in \D,$ since
		$$\sigma(\varphi_1(z),\varphi_2(z))\subseteq \sigma(\varphi_1(z))\times \sigma(\varphi_2(z))$$
		and
		$$  \varphi_1(z)\varphi_2(z)=zI_{\E_1\oplus \E_2}=\varphi_2(z)\varphi_1(z),$$
		by spectral mapping theorem we get,
		\begin{align*}
			&\sigma(\varphi_1(z),\varphi_2(z))\\
			&\subseteq \{(z_1,z_2)\in \D^2: {z_1}^{2}=\overline{\alpha}z, {z_2}^{2}=\beta z, z_1z_2=z, \text{ for  some }  \alpha,\beta \in \sigma(U_1U_2)\}.\\
			&= \{(z_1,z_2)\in \D^2: {z_1}^{2}=\overline{\alpha}z, {z_2}^{2}=\alpha z, z_1z_2=z, \text{ for  some }  \alpha \in \sigma(U_1U_2)\}.
		\end{align*}
		
		Let $z\neq 0$ and $z_1,z_2\in \D$ be such that ${z_1}^{2}=\overline{\alpha}z,\  {z_2}^{2}=\alpha z$ and $z_1z_2=z$ for some $\alpha \in \sigma (U_1U_2)$. We shall show that $\ran(\varphi_1(z)- z_1 I) +\ran(\varphi_2(z)- z_2 I) \neq \E_1\oplus \E_2$. For this, let $y\notin \ran(U_1U_2-\alpha I)$. Suppose $\begin{pmatrix}
			y \\
			0
		\end{pmatrix}\in \ran (\varphi_1(z)-z_1 I)+ \ran (\varphi_2(z)-z_2 I).$ Then there exists $\begin{pmatrix}
			h_1 \\
			h_2
		\end{pmatrix}\in \ran(\varphi_1(z)- z_1 I)$ and $\begin{pmatrix}
			k_1 \\
			k_2
		\end{pmatrix}\in \ran(\varphi_2(z)- z_2 I)$ such that $\begin{pmatrix}
			h_1 \\
			h_2
		\end{pmatrix} +\begin{pmatrix}
			k_1 \\
			k_2
		\end{pmatrix}=\begin{pmatrix}
			y \\
			0
		\end{pmatrix}.$ Since $\ran T=\ran T^*$ for any bounded normal operator $T,$ note that from \eqref{eq:ran:phi1} and \eqref{eq:ran:phi2},  we have $z_1h_1+ {U_2}^{*} h_2, z_2k_1+U_1 k_2 \in \ran (U_1U_2-\alpha I).$ Since $y\notin \ran (U_1U_2-\alpha I),$ we have $(z_1 U_1- z_2 {U_2}^{*})(k_2)\notin \ran (U_1U_2-\alpha I).$ So ${(U_1U_2-\alpha I)}^{*}(U_1(k_2))\notin \ran (U_1U_2-\alpha I).$ Which is a contradiction, as $\ran (U_1U_2-\alpha I)=\ran (U_1U_2-\alpha I)^*.$ So $\ran(\varphi_1(z)- z_1 I) +\ran(\varphi_2(z)- z_2 I) \neq \E_1\oplus \E_2$ for all $z\ne 0,$  $z_1,z_2\in \D$ such that ${z_1}^{2}=\overline{\alpha}z,\  {z_2}^{2}=\alpha z$ and $z_1z_2=z$ for some $\alpha \in \sigma (U_1U_2)$. Also, note that  $\ran(\varphi_1(0)) +\ran(\varphi_2(0)) = \E_1\oplus 0.$
		Hence, we have
		\begin{align*}
			&\sigma(\varphi_1(z),\varphi_2(z))\\&=\{(z_1,z_2)\in \D^2: {z_1}^{2}=\overline{\alpha}z, {z_2}^{2}=\alpha z, z_1z_2=z, \text{ for  some }  \alpha \in \sigma(U_1U_2)\}\\
			&=\{\pm \sqrt{z} (e^{-i \frac{\theta}{2}} , e^{i \frac{\theta}{2}}) : e^{i \theta} \in \sigma(U_1U_2) \}
		\end{align*}
		for $z\in \D.$
		
		As we saw, for any point $(z_1,z_2)\in \sigma(\varphi_1(z),\varphi_2(z))$, $z\in \D$, $$\ran(\varphi_1(z)-z_1I)+\ \ran(\varphi_2(z)-z_2I)\neq \E_1\oplus \E_2.$$ Hence $$\ran(M_{\varphi_1}-z_1I)+\ \ran(M_{\varphi_2}-z_2I)\neq H^{2}(\E_1\oplus \E_2).$$ Therefore $(z_1,z_2)\in \sigma(M_{\varphi_1},M_{\varphi_2}),$
		which implies ${\cup_{z\in \D}\sigma(\varphi_1(z),\varphi_2(z))}\subseteq\sigma(M_{\varphi_1},M_{\varphi_2}).$ Notice that for every  point $(z_1,z_2)$ in $ {\cup_{z\in \D}\sigma(\varphi_1(z),\varphi_2(z))},$ the non-singularity of $K(M_{\varphi_1}-z_1I,M_{\varphi_2}-z_2I)$ breaks at stage 3.
		
		Take $z\in \D$ and $e^{i\theta}\in \sigma(U_1U_2).$ Let $w=ze^{\frac{i \theta}{2}}, z_1=we^{-\frac{i \theta}{2}}$ and $z_2=we^{\frac{i \theta}{2}}.$ So $(z_1, z_2)\in \sigma(\varphi_1(w^{2}),\varphi_2(w^{2}))$ and $(z_1,z_2)=z(1,e^{i\theta})$. Hence we have the equality:
		\begin{align}
			&{\cup_{z\in \D}\sigma(\varphi_1(z),\varphi_2(z))}\nonumber\\
			&=\cup_{z\in \D}\{(z_1,z_2)\in \D^2: {z_1}^{2}=\overline{\alpha}z, {z_2}^{2}=\alpha z, z_1z_2=z, \text{ for some } \alpha \in \sigma(U_1U_2)\}\nonumber\\
			&=\{z(1,\alpha): z\in \D, \alpha \in \sigma (U_1 U_2)\}.\label{union}
		\end{align}
		Now as in the proof of Theorem \ref{thm:OffDiagStruc}, we have $(M_{\varphi_1}, M_{\varphi_2})$ is jointly unitarily equivalent to $(M_{z}\otimes I_{\E_2}, M_{z}\otimes U_2U_1),$ hence from \eqref{proto off} and \eqref{union}, we have
		\begin{align*}
			\sigma(M_{\varphi_1}, M_{\varphi_2})=\sigma(M_{z}\otimes I_{\E_2}, M_{z}\otimes U_2U_1)&=\{z(1,\alpha):z\in \Dc, \alpha \in \sigma (U_1U_2)\}\\
			&=\overline{\cup_{z\in \D}\sigma(\varphi_1(z),\varphi_2(z))}.
		\end{align*}
	\end{proof}
	Note that  $\sigma(V_1,V_2)=\sigma(M_{\varphi_1},M_{\varphi_2})\cup \sigma(V_1|_{\H_u},V_2|_{\H_u}),$ by Theorem \ref{thm:BCL}. Hence by Theorem \ref{offspecmulti},
	\begin{equation}\label{eq:unioninc}
		\sigma(M_{\varphi_1},M_{\varphi_2})=\overline{\cup_{z\in \D}\sigma(\varphi_1(z),\varphi_2(z))}\subseteq \sigma(V_1,V_2).
	\end{equation}
	The above inclusion is an equality if and only if $\sigma(V_1|_{\H_u},V_2|_{\H_u})\subseteq \sigma(M_{\varphi_1},M_{\varphi_2}).$

	\subsection{Range of one isometry is strictly contained in the range of other}
		
	If $(V_1,V_3)$ is  any pair of commuting isometries with $V_3$ is not unitary, set $V_2=V_1V_3.$ Then $(V_1,V_2)$ satisfy $ \ran V_2\subsetneq \ran V_1.$ Since the study of the  case $\ran V_1\subsetneq \ran V_2$ is equivalent to the case of $\ran V_2\subsetneq \ran V_1,$ we consider only the case $\ran V_2\subsetneq \ran V_1.$

	\begin{lemma}\label{lem:rangeincl} Let $(V_1,V_2)$ be a pair of commuting isometries on a Hilbert space $\H$. Then the following are equivalent:
		\begin{itemize}	
			\item[(a)] $\ran V_2\subsetneq\ran V_1.$ 		
			\item[(b)] $V_2(\ker V_1^*)\subsetneq V_1(\ker V_2^*).$ 	
			\item[(c)]	 The fringe operator $F_1=0$ and $F_2\ne 0.$
			\item[(d)] If $(\E,P,U)$ is the BCL triple for $(V_1,V_2),$ then $U(\ran P)\subsetneq\ran P^\perp$ (or equivalently $U^*(\ran P)\subsetneq\ran P^\perp$).	
			\item[(e)] If $\F$ is the defect space of $V_2$ and if $(M_z\otimes I_\F)\oplus W_2$ is the Wold decomposition of $V_2,$  then $V_1=M_\varphi\oplus W_1,$ with $\varphi(z)=W^*(Q^\perp +zQ)$ for some unitary $W$ and projection $Q\ne I$ in $\B(\F),$ and a unitary $W_1$ on $\H\ominus H^2_\D(\mathcal{F}).$	
	\end{itemize}
\end{lemma}
	
	\begin{proof} 	$(a)\iff (b)$ and $(a)\iff (c)$ follows from \eqref{eq:directsum}.\\  
		$(a)\iff (d)$: Let $(\E,P,U)$ be the BCL triple for $(V_1,V_2).$ Note that  $\ker V_1^*\subsetneq \ker V_2^*$ if and only if  $\ker M_{\phi_1}^{*}\subsetneq \ker M_{\phi_2}^{*}$ if and only if $U^{*}PU\le P^{\perp}$ and $U^{*}PU\neq  P^{\perp}$ if and only if $U(\ran P)\subsetneq \ran P^\perp,$ because $\ker M_{\phi_1}^{*}=1\otimes \ran(U^{*}PU)$ and $\ker M_{\phi_2}^{*}=1\otimes \ran P^\perp.$\\
		$(a)\iff (e)$: 	By Douglas Lemma (\cite{Douglas}), we have  $V_2=V_1V_3$ for some isometry (non-unitary) $V_3.$ Notice that $V_3$ commutes with $V_1.$ Let $(\mathcal{ F}, Q, W)$ be the BCL triple for $(V_1,V_3).$ Then, (see Theorem \ref{thm:BCL})
		$$V_1=M_{\varphi}\oplus W_1,\quad
		V_3=M_{\psi}\oplus W_3 \text{ in }\B(H^2_\D(\F)\oplus(\H\ominus H^2_\D(\mathcal{F}))),$$
		where $\phi(z)=W^*(Q^\perp +zQ)$ and $\psi(z)=(Q+zQ^\perp)W$ for all $z\in \D.$
		Now since $V_2=V_1V_3$ and $\varphi(z)\psi(z)=z$ for all $z\in \D$  we have
		$$V_2=M_z\otimes I_\F\oplus  W_2,$$
		where $W_2=W_1W_3.$  This completes the proof of $(a)\implies (e),$ because $V_3$ is not a unitary implies that $Q\ne I.$ $(e)\implies (a)$ is trivial.		
	\end{proof}

	\subsubsection{Joint spectrum}
	Consider the pair
	\begin{equation}\label{eq:Burdak pair}
		(V_1,V_2)=(U^kV^n,U^lV^m)
	\end{equation} with $V$ an isometry (not unitary) and $U$ a unitary which commutes with $V,$ for some non-negative integers $n,m,l,k$ and $n<m,$ as in \cite[Sec. 4]{Burdak}. Since $U$ is a unitary commuting with $V,$ we have $$\ran V_2=\ran U^lV^m=\ran V^m\subsetneq \ran V^n =\ran U^kV^n=\ran V_1.$$ In \cite{Burdak}, Burdak has given the model for such pairs, viz., \begin{equation}\label{model}
		V_1=(M_z^n\otimes W^k)\oplus W_1, V_2=(M_z^m\otimes W^l)\oplus W_2,
	\end{equation} for some unitary $W,W_1,W_2.$ We shall now compute the joint spectrum $\sigma(M_z^n\otimes W^k,M_z^m\otimes W^l)$ of the pure part, as an easy application of the polynomial spectral mapping theorem, viz., consider the polynomial $f:{\mathbb{C}}^{2}\rightarrow {\mathbb{C}}^{2}$ given by,$$f(z_1,z_2)=(z_1^nz_2^k, z_1^mz_2^l).$$ By the polynomial spectral mapping theorem,
		\begin{align*}\sigma(M_z^n\otimes W^k,M_z^m\otimes W^l)&=\sigma( f(M_{z}\otimes I, I\otimes W))\\
			&=f(\sigma(M_{z}\otimes I, I\otimes W) \\
			&=f(\Dc\times \sigma(W))\\
			&=\{(z^n\alpha^k,z^m\alpha^l): z\in \Dc, \alpha\in \sigma(W)\}.
		\end{align*}
Following example shows that the class considered in \cite{Burdak}, namely, the class of pairs of the type given in \eqref{eq:Burdak pair}, is only a subclass of this case.
	\begin{example}
		Let $\H=H^2_{\D^2}$ and $(V_1,V_2)=(M_{z_1},M_{z_1z_2}).$ Then there is no unitary $U$ commuting with both $V_1$ and $V_2$ such that $V_1^m=UV_2^n$ for any $m,n.$ In particular, $(V_1,V_2)$ is not of the form given in \eqref{eq:Burdak pair}. But clearly $\ran V_2\subsetneq \ran V_1.$
	\end{example}
	
	\begin{lemma}
			Let $(V_1,V_2)$ be a pair of commuting isometries, with $\ran V_2\subsetneq\ran V_1.$ Then,
		\begin{equation}\label{eq:sp:incln}
			\sigma(V_1,V_2)\subseteq \{(z_1,z_2): |z_2|\le |z_1|, z_1\in \Dc\}.
		\end{equation}
	\end{lemma}
	
	\begin{proof}
		By Douglas Lemma (\cite{Douglas}), we have  $V_2=V_1V_3$ for some isometry $V_3,$ which commutes with $V_1.$ Consider the polynomial $p:\C^2\to \C^2$ given by
		$p(z_1,z_2)=(z_1,z_1z_2).$ By the spectral mapping theorem,
		\begin{align*}
			\sigma(V_1,V_2)=\sigma(p(V_1,V_3))&=p(\sigma(V_1,V_3))\\ &\subseteq p(\Dc\times \Dc)=\{(z_1,z_1z_2): z_1,z_2\in \Dc\}\\& =\{(z_1,z_2): |z_2|\le |z_1|, z_1\in \Dc\}.
	\end{align*}	\end{proof}
	
	The inclusion in \eqref{eq:sp:incln} is sharp; see Example \ref{eg:sharp}, and it can be a strict inclusion; see Example \ref{eg:strict}.
	
	\begin{example}\label{eg:sharp}
		Consider the pair $(M_{z_1},M_{z_1z_2})$  of commuting  isometries in $H^2(\D^2).$  Notice that
		\begin{align*}
			\ran (M_{z_1z_2})&=\left\{\sum_{m,n\ge 1}^\infty a_{m,n}z_1^mz_2^n:\sum_{m,n\ge 1 }^\infty |a_{m,n}|^2<\infty \right\}\\&\subsetneq \left\{\sum_{m\ge 1,n\ge 0}^\infty a_{m,n}z_1^mz_2^n:\sum_{m\ge 1,n\ge 0 }^\infty |a_{m,n}|^2<\infty \right\}=\ran M_{z_1}.
		\end{align*}
		Consider the polynomial $p:\C^2\to \C^2$ given by
		$p(z_1,z_2)=(z_1,z_1z_2).$ By the spectral mapping theorem, we have
		\begin{align*}
			\sigma(M_{z_1},M_{z_1z_2})
		=\sigma(p(M_{z_1},M_{z_2}))&=p(\sigma(M_{z_1},M_{z_2}))= p(\Dc\times \Dc)\\
		&=\{(z_1,z_2): |z_2|\le |z_1|, z_1\in \Dc\}.
	\end{align*}
	\end{example}
		The measure of the above spectrum is non-zero, indeed it is $\frac{\pi^2}{2}.$
	
	\begin{example}\label{eg:strict}
		Let $0\le m<n.$ Consider $(M_z^m,M_z^n).$ Then, $$\ran M_z^n=\{\sum_{k=n}^\infty a_kz^k:\sum_{k=n}^\infty |a_k|^2<\infty \}\subsetneq \{\sum_{k=m}^\infty a_kz^k:\sum_{k=m}^\infty |a_k|^2<\infty \}=\ran M_z^m.$$
		By the spectral mapping theorem,
		$$\sigma(M_z^m,M_z^n)=\{(z^m,z^n): z\in \Dc\}\subsetneq \{(z_1,z_1z_2): z_1,z_2\in \Dc\}.$$
	\end{example}

\begin{example}
	Let $\H=H^2_{\D^2}\oplus H^2_{\D^2},$ and $V_i=M_{z_i}\oplus \tau_i$ for $i=1,2.$ Then $(V_1,V_2)$ is a pair of commuting isometries with defect $E_0\oplus -P_{\spn\{z_2\}}.$  Clearly the defect  is a difference of two mutually orthogonal projections and $(V_1,V_2)$  lies in the unknown case given in Table \ref{classification:table}.
\end{example}
The following example gives an irreducible  pair of commuting isometries lying in the unknown case given in Table \ref{classification:table}.
\begin{example}
	 Consider the pure pair $(V_1,V_2)$ with the BCL triple $(l^2(\Z), p_{01},\omega)$ where $\omega$ is the bilateral shift and $p_{01}$ is the projection in $l^2(\Z)$ onto $\spn\{e_0,e_1\}.$ Then
	 \begin{equation}\label{eq:difofpro}
	 	C(V_1,V_2)=E_0\otimes (\omega^*p_{01}\omega -p_{01})=E_0\otimes p_{\spn \{e_{-1}\}}-E_0\otimes p_{\spn \{e_1\}},
 	\end{equation} where $E_0$ is the one dimensional projection onto the space of constant functions in $H^2_\D$. The pair $(V_1,V_2)$ is irreducible. To show this, we shall show that the pair  $(p_{01},\omega)$  is irreducible; see \cite{Sandip,GG}. Suppose $\E_0\ne \{0\}$ is a reducing subspace for $(p_{01},\omega).$ Let $\sum_{n\in \Z}a_ne_n\in \E_0$ and $a_{n_0}\ne 0$  for some $n_0\in \Z.$ Consider
	 $${\omega^*}^{n_0}\sum_{n\in \Z}a_ne_n=\sum_{n\in \Z}a_ne_{n-n_0}\in \E_0.$$
	 This implies that $p_{01}\sum_{n\in \Z}a_ne_{n-n_0}=a_{n_0}e_0+a_{n_0+1}e_1\in \E_0.$ Hence $\omega^*(a_{n_0}e_0+a_{n_0+1}e_1)=a_{n_0}e_{-1}+a_{n_0+1}e_0\in \E_0.$ So $a_{n_0}e_{-1}\in \E_0.$ Thus $e_{-1}\in \E_0,$ shows that $\E_0=l^2(\Z).$
	
	 The pair $(V_1,V_2)$ lies in the unknown case of $\H_i\ne 0$ for all $i=1,2,3,4$ mentioned in Table \ref{classification:table}. To see this first note that \eqref{eq:difofpro} shows that $C(V_1,V_2)$ is a difference of two mutually orthogonal projections and one can see that, it is not in any other case of the Table \ref{classification:table}, using the characterization in terms of BCL given for those cases.
\end{example}

\textsf{Acknowledgement}: The authors are thankful to the referee for a careful reading and pertinent comments. The first author is supported by the J C Bose Fellowship JCB/2021/000041 of SERB. The third author gratefully acknowledges support from the Indian Institute of Science through the Institute of Eminence post-doctoral fellowship and from the University Grants Commission, India through a D. S. Kothari post-doctoral fellowship. We are thankful to Dan Popovici for bringing a good amount of existing literature to our notice.

\textbf{Data availability statement}: Data sharing is not applicable to this article as no data sets were generated or analysed during the current study.
	

\end{document}